\documentclass[11pt]{article}

\usepackage{amsmath,amsfonts,amssymb}
\usepackage{mathdots}
\usepackage{graphicx}
\usepackage[colorlinks,linkcolor=blue,citecolor=blue]{hyperref}
\usepackage{color}

\newtheorem{theorem}{Theorem}[section]
\newtheorem{proposition}{Proposition}[section]
\newtheorem{lemma}[theorem]{Lemma}
\renewcommand{\Vec}[1]{\textup{vec}\left(#1\right)}

\begin{document}

\title{Solutions of the matrix equation $p(X)=A$,\\ with polynomial function $p(\lambda)$ \\ over field extensions of~$\mathbb{Q}$}

\author{G.J.~Groenewald\footnotemark[1]
, D.B.~Janse van Rensburg\footnotemark[1],
A.C.M.~Ran\footnotemark[2],\\
F.~Theron\footnotemark[1],
M.~van~Straaten\footnotemark[1] \footnotemark[3]}
\renewcommand{\thefootnote}{\fnsymbol{footnote}}
\footnotetext[1]{School~of~Computer,~Statistical~and~Mathematical~Sciences,
North-West~University,
Research Focus: Pure and Applied Analytics,
Private~Bag~X6001,
Potchefstroom~2520,
South Africa.
E-mail: \texttt{
gilbert.groenewald@nwu.ac.za,  dawie.jansevanrensburg@nwu.ac.za,frieda.theron@nwu.ac.za, madelein.vanstraaten@nwu.ac.za}
}
\footnotetext[2]{Department of Mathematics,Vrije Universiteit Amsterdam, De Boelelaan
    1111, 1081 HV Amsterdam, The Netherlands
    and Research Focus: Pure and Applied Analytics, North-West~University,
Potchefstroom,
South Africa. E-mail:
    \texttt{a.c.m.ran@vu.nl}}
\footnotetext[3]{DSI-NRF Centre of Excellence in Mathematical and Statistical Sciences (CoE-MaSS)}    
\date{}

\maketitle

\noindent
{\bf Abstract.}

\noindent
Let $\mathbb{H}$ be a field with $\mathbb{Q}\subset\mathbb{H}\subset\mathbb{C}$, and
let $p(\lambda)$ be a polynomial in $\mathbb{H}[\lambda]$, and let $A\in\mathbb{H}^{n\times n}$  be nonderogatory. In this paper we consider the problem of finding a solution $X\in\mathbb{H}^{n\times n}$ to $p(X)=A$. A necessary condition for this to be possible is already known from \cite{Drazin}. Under an additional condition we provide an explicit construction of such solutions.
The similarities and differences with the derogatory case will be discussed as well.

One of the tools needed in the paper is a new canonical form, which may be of independent interest. It combines elements of the rational canonical form with elements of the Jordan canonical form.

\bigskip
\noindent
\emph{Kewords:}
solutions of polynomial matrix equations, canonical forms, companion matrices, matrices over field extensions of the rationals, linear matrix equations, nonderogatory matrices


\section{Introduction}

Let  $\mathbb{H}$ be a field such that $\mathbb{Q}\subset\mathbb{H}\subset\mathbb{C}$. In this paper we consider the equation $p(X)=A$, where $p(\lambda)$ is a polynomial of degree $\ell$ with coefficients in $\mathbb{H}$ and $A$ is an $n\times n$ matrix with entries in $\mathbb{H}$. The problem we are considering is whether or not a solution $X$ with entries in $\mathbb{H}$ exists and if it does, how to construct it. The paper is highly motivated by \cite{Drazin}. In that paper the case where $A$ is a simple matrix is solved. The main result of our paper is to extend these results to the case where $A$ is nonderogatory, and it can be stated in a somewhat imprecise way as follows. (Recall that a matrix $A$ is nonderogatory if for each eigenvalue the geometric multiplicity of that eigenvalue is one.) The first part of the following theorem already appears in \cite{Drazin}.

\begin{theorem}\label{thm:main_intro}
Let $A$ be an $n\times n$ matrix with entries in a field $\mathbb{H}$, with $\mathbb{Q}\subset \mathbb{H}\subset\mathbb{C}$, and let $p_A(\lambda)=g_1(\lambda)^{d_1}\cdot g_2(\lambda)^{d_2}\cdots g_r(\lambda)^{d_r}$ be the factorization of the characteristic polynomial of $A$ with $g_j(\lambda)$'s pairwise coprime monic and irreducible polynomials. In particular, $A$ is nonderogatory.
Let $p(\lambda)\in\mathbb{H}[\lambda]$. 
If $p(X)=A$ has a solution, then for each root $\lambda_j$ of some $g_k(\lambda)$ there is a solution $\mu\in \mathbb{H}(\lambda_j)$ with $p(\mu)=\lambda_j$. 

Conversely, if $p(\mu)=\lambda_j$ has a solution in $\mathbb{H}(\lambda_j)$ for each eigenvalue $\lambda_j$ of $A$, then there is a solution $X$ of $p(X)=A$ with entries in $\mathbb{H}$ provided a certain finite set of linear matrix equations is solvable, and a solution can be  constructed explicitly.
\end{theorem}

To illustrate the theorem, in particular the second part, consider the following example.

\noindent
{\bf Example.} Let $A=\begin{bmatrix} 4 & 1 \\ 0 & 4 \end{bmatrix}$ and $p(\lambda)=\lambda^2$, with $\mathbb{H}=\mathbb{Q}$.  By a direct calculation one can see that $X^2=A$ has the solutions  $X=\pm \begin{bmatrix} 2 & \tfrac{1}{4}\\ 0 & 2\end{bmatrix}$. 
Our methods will lead to the fact that if a solution exists, it would have to be an upper triangular Toeplitz matrix $X=\begin{bmatrix} x_1 & x_2 \\ 0 & x_1\end{bmatrix}$. It is then immediate that $x_1$ should be a solution to $p(x_1)=4$, so that $x_1=\pm 2$, and that $x_2$ should satisfy the linear equation $x_1x_2+x_2x_1=1$. {\phantom{.}}\hfill $\Box$

\medskip

The theorem will be made more precise and proved in Section \ref{sec:solving} of the paper, see Theorem \ref{thm:main}. The algorithm for finding the solution $X$ of $p(X)=A$ mentioned in the theorem will be described in detail in that section. We shall also extend it partly to the case of a general derogatory matrix $A$. For the latter case we  present sufficient conditions for the existence of a solution $X\in\mathbb{H}^{n\times n}$ to $p(X)=A$. 

Preceding the proof of the theorem we present in several sections different tools that are used in the proof. In particular, in Section 2 we develop a canonical form which is new to the best of our knowledge, and which is, we think, of independent interest.  We view this canonical form as the second main result of the paper.

In order to describe this canonical form, let $C_g$ denote the $k\times k$ companion matrix corresponding to a polynomial $g(\lambda)$ of degree $k$. For an $n\times n$ matrix $A$ with entries in a field $\mathbb{H}$ consider its characteristic polynomial $p_A(\lambda)$, and its factorization into invariant factors
$p_A(\lambda) =\prod_{j=1}^l g_j(\lambda)^{d_j}$ with each $g_j(\lambda)$ an irreducible monic polynomial of degree $k_j$. Then $A = T^{-1}CT$, where $T$ is an invertible matrix with entries in $\mathbb{H}$, and $C$ is a direct sum of matrices $C_1, \ldots , C_l$, where each $C_j$ is a $d_jk_j\times  d_jk_j$ block upper triangular matrix of the form $C_j=I_{d_j}\otimes C_{g_j} +N\otimes I_{k_j} $, where $ N$ is the $d_j\times d_j$ upper triangular matrix with zeroes everywhere except for ones in the $(j,j+1)$ entries. See Theorem \ref{thm:companion-Jordanform} below.

Let $X$ be a  block upper triangular matrix, that is, $X=\begin{bmatrix} W & Z\\ 0 & Y \end{bmatrix}$.
For a polynomial $p(\lambda)$, 
the form of $p(X)$ is given by
$$
p(X)=\begin{bmatrix} p(W) & \tilde{Z} \\ 0 & p(Y)\end{bmatrix}
$$
for some matrix $\tilde{Z}$ which depends on $W, Y,Z$ and the polynomial $p(\lambda)$ in an intricate manner.
Introducing the notation $\tilde{Z}=\Delta p(W,Y)(Z)$ it turns out that this has many properties reminiscent of differentiation (see e.g., \cite{kv-v}). This is further discussed in Section 3, and will play a prominent role in later sections.

In Section 4 we recall several well known results on linear matrix equations which will play a role in the sequel.

In Section 5 the discussion on solutions of $p(X)=A$ for a nonderogatory matrix $A$ is started. Several properties of the solutions are presented. In particular, it is shown that if $A$ is written in companion-Jordan form as in Section 2, then $X$ must have a compatible decomposition. More precisely, if $A=T^{-1}(C_1\oplus \cdots \oplus C_r)T$ is the companion-Jordan form of $A$, then $X=T^{-1}(X_1\oplus \cdots \oplus X_r)T$ with $X_j$ commuting with $C_j$. This forces $X_j$ to have the form of a block upper triangular Toeplitz matrix. 

At this point in the discussion we present  in Section 6 a review of the main results of \cite{Drazin}. It is shown in \cite{Drazin} that if $A$ is nonderogatory with entries in a field $\mathbb{H}$, and $p(\lambda)\in\mathbb{H}[\lambda]$, then the following holds. If $p(X)=A$ has a solution $X$ with entries in $\mathbb{H}$, then for every eigenvalue $\lambda_0$ of $A$ the equation $p(\mu)=\lambda_0$ must have a solution in the field $\mathbb{H}(\lambda_0)$. Conversely, under the extra assumption that $A$ is a simple matrix, if $p(\mu)=\lambda_0$ has a solution in $\mathbb{H}(\lambda_0)$ for every eigenvalue $\lambda_0$ of $A$, then $p(X)=A$ has a solution $X$ with entries in $\mathbb{H}$, which can be constructed explicitly.  

Based on the construction in \cite{Drazin} for the case of a simple matrix, we extend that construction, using the tools developed in the earlier sections, to the nonderogatory case in Section 7. This leads us to the main theorem of the paper as presented above. The (im)possibility to extend the results to the derogatory case is discussed as well.

Finally, in Section 8 we present the special case where $p(\lambda)=\lambda^\ell$ and in Section 9 we present several illustrative examples.

As far as we are aware, there is only a limited literature on the problem of polynomial equations of the form considered here. In fact, there is a substantial number of papers which are concerned with the special case of the $m$th root. For a selection of relatively recent papers on this problem, see \cite{otero,psarrakos,reams, tenhave} and the sources quoted therein; for early sources we refer to \cite{roth, wedderburn} and the sources quoted therein. In contrast, we know of only a few papers that are specifically concerned with the problem of solving $p(X)=A$ where $p(\lambda)$ is a polynomial and $X$ and $A$ are matrices. One of the earliest sources is the paper by W.E. Roth \cite{roth}. The focus there is on finding solutions $X$ which are polynomials in $A$. The underlying field is not explicitly stated, but from the text it is clear that this is supposed to be the complex field. The paper is worth reading also for its introduction, which provides references to the very early literature in this area, going back to papers of Cayley and Sylvester dating from the second half of the nineteenth century. In much more modern terms are the papers by E. Spiegel \cite{Spiegel} and by J-C. Evard and F. Uhlig \cite{EU}. Evard and Uhlig deal with the more general situation $f(X)=A$, where $f(\lambda)$ is a function which is a complex holomorphic function on some open domain in the complex plane, and $A$ and $X$ are complex matrices, whereas Spiegel focusses on nonconstant complex polynomials. In \cite{EU} some attention is given to the real case as well. Finally, the paper by M.P. Drazin \cite{Drazin} is, as far as we know, the first paper to consider the problem on field extensions of $\mathbb{Q}$. It is from Drazin's presentation that we pick up the story.

\section{The companion-Jordan form}\label{sec:compJordan}

Let $g (\lambda)=\lambda^n+a_{n-1}\lambda^{n-1} + \cdots +a_0$ be a monic polynomial in $\mathbb{H}[\lambda]$ and let 
$$
C_g=\begin{bmatrix} 
0 & 1 & 0 & \cdots & 0 \\
0 & 0 & 1 & \ddots & \vdots \\
\vdots & \ddots & \ddots & \ddots &0\\
0 & \cdots & \cdots & 0 & 1\\
-a_0 & -a_1 & \cdots & \cdots & -a_{n-1}
\end{bmatrix}
$$
 be its $n\times n$ companion matrix.

\begin{theorem}{\rm (Companion-Jordan form.)}\label{thm:companion-Jordanform} 
Let $A\in\mathbb{H}^{n\times n}$ and let $$
p_A(\lambda)= g_1(\lambda)^{d_1} \cdot g_2(\lambda)^{d_2} \cdots g_r(\lambda)^{d_r}
$$ be the factorization of its characteristic polynomial into elementary factors, with each $g_j(\lambda)$ an irreducible monic  polynomial of degree $k_j$, so $\sum_{j=1}^r k_jd_j=n$. Then there is an invertible matrix $T\in \mathbb{H}^{n\times n}$ such that 
\begin{equation}\label{eq:cJ0}
A=T^{-1} (C_1\oplus C_2\oplus \cdots \oplus C_r)T,
\end{equation}
where each $C_j$ is a $d_jk_j\times d_jk_j$ block upper triangular matrix of the form
\begin{equation}\label{eq:companion-Jordanform}
C_j=\begin{bmatrix} 
C_{g_j} & I_{k_j} & 0 & \cdots & 0\\
0 & C_{g_j} & I_{k_j} & \ddots & \vdots \\
\vdots & \ddots & \ddots & \ddots &0\\
\vdots  &  & \ddots & C_{g_j}& I_{k_j}\\
0 & \cdots & \cdots & 0 & C_{g_j}
\end{bmatrix}.
\end{equation}
The form \eqref{eq:cJ0} is unique up to permutation of the blocks $C_1, \ldots , C_r$ in the direct sum.
\end{theorem}

The matrix  $C_1\oplus C_2 \oplus \cdots \oplus C_r$ in \eqref{eq:cJ0}, with each $C_j$ as in equation \eqref{eq:companion-Jordanform} will be called the \emph{companion-Jordan form} of $A$. Note that the polynomials $g_j(\lambda)$ need not be distinct, in fact, if $A$ has several Jordan blocks in the Jordan canonical form corresponding to the same eigenvalue, then there will be as many $g_j(\lambda)$ with that eigenvalue as a root as there are Jordan blocks. The remainder of this section will be devoted to the proof of Theorem~\ref{thm:companion-Jordanform}.

We first prove the assertion of the theorem for $p_{A}(\lambda)$ a power of a single irreducible polynomial, say $p_{A}(\lambda) = g(\lambda)^{d}$, where $g(\lambda) \in \mathbb{H}[\lambda]$ is monic and irreducible. The companion matrix of $p_{A}(\lambda)$ is denoted by $C_{g^{d}}$ and for brevity we will use $q(\lambda)$ instead of $p_{A}(\lambda)$. We show that there is an invertible matrix $T\in \mathbb{H}^{dn\times dn}$ such that
$$
T^{-1}C_{g^d}T=
\begin{bmatrix}
C_g & I_n & 0 & \cdots & 0\\
0 & C_g & I_n & \ddots & \vdots \\
\vdots & \ddots & \ddots & \ddots &0\\
\vdots  & \cdots & \ddots & C_g& I_n\\
0 & \cdots & \cdots & 0 & C_g
\end{bmatrix}.
$$
 For the construction of the matrix $T$ we need an intermezzo on the Jordan chains of companion matrices.

\subsubsection*{Jordan chains of companion matrices}

For any $\lambda$  and any $n>1$, let $v_n(\lambda)$ be the vector
$$
v_n(\lambda)=\begin{bmatrix} 1 \\ \lambda \\ \lambda^2 \\ \vdots \\ \lambda^{n-1}\end{bmatrix}.$$
It is well known and easily checked that if $\lambda_0$ is a root of $g(\lambda)$, so $g(\lambda_0)=0$, and $g(\lambda)$ has degree $n$, then 
$C_g v_n(\lambda_0)=\lambda_0 v_n(\lambda_0)$, and in fact $v_n(\lambda_0)$ spans the eigenspace of $C_g$ corresponding to $\lambda_0$. In particular, $C_g$ is nonderogatory.

Consider $C_q v_n(\lambda)$ where $q(\lambda)=g(\lambda)^d$ has degree $n$ and view this as a vector function $h(\lambda)$ of $\lambda$.
Then
$$
h(\lambda)= C_qv_n(\lambda)=\begin{bmatrix} \lambda \\ \lambda^2 \\ \vdots \\ \lambda^{n-1} \\ -q(\lambda)+\lambda^n\end{bmatrix}=\lambda v_n(\lambda) -q(\lambda) e_n.
$$
So, the derivatives $h^{(k)}(\lambda)$ for $k=1,2, \ldots $ are given by
$$
h^\prime (\lambda) =C_qv_n^\prime(\lambda)= v_n(\lambda)+\lambda v_n^\prime(\lambda) -q^\prime(\lambda) e_n
$$
and
$$
h^{\prime\prime}(\lambda)=C_qv_n^{\prime\prime}(\lambda)=
2v_n^\prime(\lambda) + \lambda v_n^{\prime\prime}(\lambda) -q^{\prime\prime}(\lambda)e_n,
$$
while in general, one proves by induction that
$$
h^{(k)}(\lambda)=C_qv_n^{(k)}(\lambda)=k v_n^{(k-1)}(\lambda) +\lambda v_n^{(k)}(\lambda) -q^{(k)}(\lambda) e_n.
$$

Now suppose that $\lambda_0$ is a root of $g(\lambda)$, then it is a root of $q(\lambda)$ of multiplicity $d$, and hence $q^{(k)}(\lambda_0)=0$ for $k=0, 1, \ldots , d-1$. From the above formulas we obtain that
\begin{equation}\label{eq:Jordanchainv_n}
v_n(\lambda_0), \ v_n^\prime(\lambda_0),\  \frac{1}{2}v_n^{\prime\prime}(\lambda_0),\  \frac{1}{3!}v_n^{(3)}(\lambda_0),\ \ldots\ ,\ \frac{1}{(d-1)!}v_n^{(d-1)}(\lambda_0)
\end{equation}
is a Jordan chain of length $d$ of $C_q$ corresponding to the eigenvalue $\lambda_0$.

We rewrite the Jordan chain of length $d$ of $C_q$ in terms of the following $n\times n$ matrices: let
\begin{equation}\label{eqS_j}
S_j =\begin{bmatrix} 0_{j\times n-j} & 0_{j\times j} \\[4mm]
{\rm diag\,}\begin{pmatrix} k \\ j\end{pmatrix}_{k=j}^{n-1} & 0_{n-j\times j}\end{bmatrix} \quad\mbox{\ for\ } j=1, \ldots , n.
\end{equation}
For instance,
$$
S_1=\begin{bmatrix} 0 & \cdots& \cdots & \cdots & 0 \\
1 & 0 &  & & \vdots \\
0 & 2 & \ddots & & \vdots \\
\vdots & \ddots  & \ddots &\ddots  & \vdots \\
0 & \cdots & 0 & n-1 & 0\end{bmatrix}.
$$
Then it is easily seen that $S_1v_n(\lambda_0)=v_n^\prime(\lambda_0)$, $S_2v_n(\lambda_0) =\frac{1}{2} v_n^{\prime\prime}(\lambda_0)$, and in general $S_jv_n(\lambda_0)=\frac{1}{j!}v_n^{(j)}(\lambda_0)$. Hence 
$$
v_n(\lambda_0), S_1v_n(\lambda_0), \ldots , S_{d-1}v_n(\lambda_0)
$$
is a Jordan chain of length $d$ of $C_q$ corresponding to the eigenvalue $\lambda_0$, so
\begin{equation}\label{eq:Jordanchainlambda0}
C_q v_n(\lambda_0)=\lambda_0 v_n(\lambda_0) \mbox{\ \  and\ \ } C_q S_jv_n(\lambda_0)= \lambda_0 S_jv_n(\lambda_0)+S_{j-1}v_n(\lambda_0),
\end{equation}
where we put $S_0=I_n$.  

\subsubsection*{Reduction to Companion-Jordan form}

We now return to the situation where  $p_A(\lambda)=q(\lambda)=g(\lambda)^d$, with $g(\lambda)=g_0+g_1 \lambda+\cdots g_{k-1}\lambda^{k-1} +\lambda^k$ a polynomial over some field $\mathbb{H}$, which is
 monic and irreducible over $\mathbb{H}$.
Let $\lambda_0$ be a complex root of $g(\lambda)$.
Then the multiplicity of $\lambda_0$ as an eigenvalue of $A$ is $d$. Note that every root of $g(\lambda)$ has the same multiplicity $d$ as a root of $q(\lambda)$.  Denote the degree of $g(\lambda)$ by $k$, so $n=k\cdot d$.

Next, we have for $j\geq k$ that $\lambda_0^j$ can be expressed as a linear combination of $1, \lambda_0, \ldots , \lambda_0^{k-1}$ with coefficients in $\mathbb{H}$, since $\lambda_0$ is a (complex) root of $g(\lambda)$. Following Drazin \cite{Drazin}
we now  construct an $(n-k)\times k$ matrix $W_1$ with entries in $\mathbb{H}$ such that 
$$
v_n(\lambda_0)= \begin{bmatrix} I_k \\ W_1\end{bmatrix} v_{k}(\lambda_0).
$$
We denote the matrix $\begin{bmatrix} I_k \\ W_1\end{bmatrix}$ by $W$. Note that $W$ is an $n\times k$ matrix, 
and that the matrix $W$ is independent of the specific root of $g(\lambda)$ which was chosen a priori. Now let $V_g=\begin{bmatrix} v_k(\lambda_1) & \cdots & v_k(\lambda_k)\end{bmatrix}$ be the  $k\times k$ Vandermonde matrix obtained from the $k$ roots $\lambda_1, \ldots , \lambda_k$ of $g$. Since irreducibility of $g$ implies that these $k$ roots are different from each other (see \cite{pinter}, Chapter 31, Theorem 1), we have that $V_g$ is invertible. Then 
$$
\begin{bmatrix} v_n(\lambda_1) & \cdots & v_n(\lambda_k)\end{bmatrix} =WV_g.
$$
Put $V_{n,g}=\begin{bmatrix} v_n(\lambda_1) & \cdots & v_n(\lambda_k)\end{bmatrix}$.

Next, we construct the matrix $T=\begin{bmatrix} W & S_1W & S_2 W & \cdots & S_{d-1} W\end{bmatrix}$. Note that $T$ is an $n\times n$ matrix with entries in $\mathbb{H}$.

By a dimension count we now have that this is an $n\times n$ matrix where $n$ is the degree of $q(\lambda)$. 
Notice that $v_n(\lambda_0)= W v_k(\lambda_0)$ by construction of $W$, and hence
$v_n^\prime(\lambda_0)=Wv_k^\prime(\lambda_0)$.

We then have
\begin{equation}\label{eq:compJordanoneblock}
C_qT=T\begin{bmatrix} 
C_g & I_n & 0 & \cdots & 0\\
0 & C_g & I_n & \ddots & \vdots \\
\vdots & \ddots & \ddots & \ddots &0\\
\vdots  &  & \ddots & C_g& I_n\\
0 & \cdots & \cdots & 0 & C_g
\end{bmatrix}.
\end{equation}
Indeed, 
$$
C_qT= \begin{bmatrix} C_qW & C_qS_1W & \cdots & C_qS_{d-1}W\end{bmatrix}.
$$
We shall show that $C_qW=WC_g$, and $C_qS_jW=S_{j-1}W+S_jWC_g$ for $j=1, \ldots , d-1$ (where $S_0=I$). Then \eqref{eq:compJordanoneblock} follows.

To show that $C_qW=WC_g$, multiply on the right with $V_g$. Since $V_g$ is invertible, it is enough to show that $C_qWV_g=WC_gV_g$. Denote the diagonal $k\times k$ matrix with eigenvalues $\lambda_1, \ldots , \lambda_k$ by $D_g$, then $C_gV_g=V_gD_g$. Hence $WC_gV_g=WV_gD_g$. 
Now $WV_g=V_{n,g}$, so $WC_gV_g=V_{n,g}D_g$, and since the columns of $V_{n,g}$ are eigenvalues of $C_q$, we have that this is equal to $C_qV_{n,g}=C_qWV_g$ as desired.

To show that $C_qS_jW=S_{j-1}W+S_jWC_g$ for $j=1, \ldots , d-1$, again multiply on the right with $V_g$. It suffices to show that $C_qS_jWV_g=S_{j-1}WV_g+S_jWC_gV_g$ for $j=1, \ldots , d-1$. 
Using $WV_g=V_{n,g}$ and $C_qWV_g=WC_gV_g=V_{n,g}D_g$ we have from \eqref{eq:Jordanchainlambda0}, applied to all roots of $g$, the following
$$
C_qS_jWV_g=C_q S_j V_{n.g} = S_jV_{n,g}D_g+S_{j-1}V_{n,g} =S_jWC_gV_g+S_{j-1}WV_g,
$$
as desired.

Finally, we turn to a matrix where the characteristic polynomial has more than one distinct irreducible factor. We use the fact that any $n\times n$ matrix $A$  with entries in a field $\mathbb{H}$ is similar to a block direct sum of companion matrices, with similarity matrix in $\mathbb{H}^{n\times n}$. This is called the rational canonical form, see e.g., Hoffman and Kunze \cite[Section~7.2, p.199]{HK} and Roman \cite[Theorems~7.14 and 7.16]{roman}. Thus Theorem~\ref{thm:companion-Jordanform} is proved except for the uniqueness statement.

Uniqueness up to permutation of the blocks can be seen by realizing that the matrix 
\eqref{eq:companion-Jordanform} is similar to the matrix
$$
\begin{bmatrix} 
D_{g_j} & I_{k_j} & 0 & \cdots & 0\\
0 & D_{g_j}& I_{k_j} & \ddots & \vdots \\
\vdots & \ddots & \ddots & \ddots &0\\
\vdots  & \cdots & \ddots & D_{g_j}& I_{k_j}\\
0 & \cdots & \cdots & 0 & D_{g_j}
\end{bmatrix},
$$
using similarity with a block diagonal matrix with $V_{g_j}$ on the diagonal. In turn, this matrix is similar by permutation to $J_{d_j}(\lambda_1) \oplus \cdots \oplus J_{d_j} (\lambda_{k_j})$, where $\lambda_1 , \ldots , \lambda_{k_j}$ are the distinct roots of $g_j(\lambda)$. Then combining the fact that the rational canonical form is unique up to permutation of its blocks (see \cite{HornJohn1}, p. 200) with the fact that the Jordan canonical form is unique up to permutation of its blocks (see \cite{HornJohn1}, p. 167), we see that the companion-Jordan form is also unique up to permutation of its blocks.

\section{Computation of the action of a polynomial $p(\lambda)$ on a block-upper triangular matrix}\label{sec:deltas}

For use in the  sections to come, we discuss here the action of a function $p$ from a ``nice enough set of functions'' on a block upper triangular matrix of the form $X=\begin{bmatrix} W & Z\\ 0 & Y \end{bmatrix}$. For a precise definition of what ``nice enough'' means here, we can refer to Kaliuzhnyi-Verbovetskyi and Vinnikov \cite{kv-v}, but analytic in a neighbourhood of the spectrum will work, and so certainly it will work for polynomials.

One obtains that $p(X)$ is block upper triangular: 
$$
p(X)=\begin{bmatrix} p(W) & \tilde{Z} \\ 0 & p(Y)\end{bmatrix},
$$
for some matrix $\tilde{Z}$. Following notation of the book by Kaliuzhnyi-Verbovetskyi and Vinnikov \cite{kv-v}, we shall denote the matrix $\tilde{Z}$ by $\Delta p(W,Y)(Z)$. Then this operator has the following properties, see for example Sections 2.3 and 2.4 in \cite{kv-v}:
\begin{itemize}
\item $\Delta p(W,Y)(Z)$ is linear in $Z$, that is, \\ $\Delta p(W,Y)(a_1Z_1+a_2Z_2)=
a_1\Delta p(W,Y)(Z_1)+a_2\Delta p(W,Y)(Z_2)$,
\item
$\Delta p(W,Y)(Z)$ is linear in $p(\lambda)$,
\item if $p(\lambda)$ is a constant times the identity, then $\Delta p(W,Y)(Z)=0$,
\item the product rule holds:\\ $\Delta (pq)(W,Y)(Z)=p(W)\Delta q(W,Y)(Z)+\Delta p(W,Y)(Z)q(Y)$,
\item the chain rule holds:\\ $\Delta (p\circ q) (W,Y)(Z)= \Delta p(q(W),q(Y))(Z) \Delta q(W,Y)(Z)$,
\item there is an ``inverse formula'': defining $p^{-1}(X)=p(X)^{-1}$ whenever this exists, then $\Delta p^{-1} (W,Y)(Z)=-p^{-1}(W) \Delta p(W,Y)(Z) p^{-1}(Y)$,
\item a divided-difference result holds:\\ $p(X_1)-p(X_2)=\Delta p(X_1,X_2)(X_1-X_2)$.
\end{itemize}
So in many ways, the properties of the operator $\Delta$ are reminiscent of those of differentiation.

We compute explicitly some examples.
We can find $\Delta q_m(W,Y)(Z)$ for the function $q_m(\lambda)=\lambda^m$. Consider first $q_m(X)$ for $m=2$, then from
$$
X^2=\begin{bmatrix} W^2 & WZ+ZY \\ 0 & Y^2\end{bmatrix}
$$
we see
$$
\Delta q_2(W,Y)(Z)=WZ+ZY.
$$
This is an invertible linear map on the set of matrices if and only if $W$ and $-Y$ have no mutual eigenvalues, see for example Theorem 4.4.6, p.\ 270 in Horn and Johnson \cite{HornJohn} and Theorem 12.3.2, p.\ 414 in Lancaster and Tismenetsky~ \cite{LT}.

For $m=3$, so $q_3(\lambda)=\lambda^3$, we have in a similar manner
$$
\Delta q_3(W,Y)(Z)=W^2Z+\Delta q_2(W,Y)(Z)\cdot Y=W^2Z+WZY+ZY^2
$$
and likewise, for $m=4$:
$$
\Delta q_4(W,Y)(Z)=W^3Z+\Delta q_3(W,Y)(Z)\cdot Y=W^3Z+W^2ZY+WZY^2+ZY^3.
$$

Continuing in this way, we see that in general we have
\begin{align*}
\Delta q_m(W,Y)(Z)&=W^{m-1}Z+\Delta q_{m-1}(W,Y)(Z)\cdot Y\\
  &= W^{m-1}Z+W^{m-2}ZY+W^{m-3}ZY^2+\cdots +ZY^{m-1}\\
  &= \sum_{j=0}^{m-1}W^{m-1-j}ZY^j,
\end{align*}
which is linear in $Z$. 

Next, consider a general polynomial $p(\lambda)=\sum_{m=0}^\ell p_m\lambda^m$, with $p_m$ in some field $\mathbb{H}$. If $X=\begin{bmatrix} W & Z\\ 0 & Y \end{bmatrix}$, then 
$$
p(X)=\begin{bmatrix}
p(W) & \Delta p(W,Y)(Z) \\
 0   &    p(Y)
\end{bmatrix},
$$
where $\Delta p(W,Y)(Z)$ is as follows
\begin{align*}
\Delta p(W,Y)(Z)&=\sum_{m=1}^\ell p_m \Delta q_m(W,Y)(Z)\\
      & =\sum_{m=1}^\ell p_m \sum_{j=0}^{m-1} W^{m-1-j}ZY^j.
\end{align*}
Observe that this is a linear function in $Z$. Note also that if $W, Y, Z$ have entries in a field $\mathbb{H}$, and $p(\lambda)\in \mathbb{H}[\lambda]$, then $\Delta p(W,Y)(Z)$ has entries in $\mathbb{H}$.

\section{Linear matrix equations}\label{sec:lineq}

Let $A=(a_{ij})\in\mathbb{H}^{m\times n}$ and $B=(b_{ij})\in\mathbb{H}^{p\times q}$. Then the Kronecker product of $A$ and $B$ is defined (see \cite{HornJohn}, Definition~4.2.1, see also Section 12.1, p.\ 407, of \cite{LT}) as
\begin{equation*}
A\otimes B=\begin{bmatrix}
a_{11}B&\cdots&a_{1n}B\\
\vdots&\ddots&\vdots \\
a_{m1}B&\cdots&a_{mn}B
\end{bmatrix}\in\mathbb{H}^{mp\times nq}.
\end{equation*}

Let $A$ be an $m\times n$ matrix over $\mathbb{H}$ with columns $a_1,\ldots,a_n$. Then the mapping $\text{vec}:\mathbb{H}^{m\times n}\rightarrow \mathbb{H}^{mn}$ is defined by
\begin{equation*}
\Vec{A}=\begin{bmatrix}
a_1^T & a_2^T & \cdots & a_n^T
\end{bmatrix} ^T.
\end{equation*}

\begin{lemma}[{\cite[Lemma~4.3.1]{HornJohn}}]\label{lemAXB=C}
Let $A\in\mathbb{H}^{m\times n}$, $B\in\mathbb{H}^{p\times q}$ and $C\in\mathbb{H}^{m\times q}$ be given and let $X\in\mathbb{H}^{n\times p}$ be unknown. The matrix equation $AXB=C$ is equivalent to the system of $qm$ linear equations in $np$ unknowns given by
$$(B^T\otimes A)\,\Vec{X}=\Vec{C}.$$
\end{lemma}

This result can also be found in \cite{LT}, Proposition 12.1.4, p.\ 410.

Later on we shall consider equations of the type $\sum_{j=1}^k a_jA_jXB_j = C$, which, using linearity of the Kronecker product, is equivalent to the following system of equations: 
\begin{equation}\label{eq:veceq}
\left(\sum_{j=1}^k a_j(B_j^T \otimes A_j) \right) \,\Vec{X}=\Vec{C}.
\end{equation}
Now suppose that the entries of the matrices $A_j, B_j, C$, as well as the constants $a_j$, are all in a field $\mathbb{H}$, and suppose that $\Vec{C}$ is contained in the column space of $\sum_{j=1}^k a_j(B_j^T \otimes A_j)$. Then there is a solution $X$ over $\mathbb{H}$. This is due to the fact that the
reduction to row-reduced echelon form of the augmented matrix $\begin{bmatrix} \sum_{j=1}^k a_j(B_j^T \otimes A_j) & \Vec{C}\end{bmatrix}$ requires only operations in $\mathbb{H}$. Then following
the standard ``Algorithm V" in \cite{Lay} we see that there exists a basis in the null space
of the coefficient matrix in the system \eqref{eq:veceq} of $qm$ linear equations in $np$ unknowns
and the existence of a solution in the field $\mathbb{H}$
can be produced as in Step 2 of Algorithm V in \cite{Lay}.

\section{Properties of solutions of $p(X)=A$ for $A$ nonderogatory}\label{sec:properties}

In this section we prove the following proposition.

\begin{proposition}\label{prop:formofX}
Consider a nonderogatory matrix $A$ with entries in a field $\mathbb{H}$. Suppose that the characteristic polynomial $p_A(\lambda)$ factors as $p_A(\lambda)=g_1(\lambda)^{d_1} \cdot g_2(\lambda)^{d_2} \cdots g_r(\lambda)^{d_r}$ with each $g_j(\lambda)$ a monic irreducible polynomial over $\mathbb{H}$, and suppose that the polynomials $g_j(\lambda)$ are coprime. Write
$$
A=T^{-1} (C_1\oplus C_2\oplus\cdots \oplus C_r)T,
$$
with 
$C_j$ as in \eqref{eq:companion-Jordanform}.
Then any solution $X$ of $p(X)=A$ must be of the form
$$
X=T^{-1} (X_1\oplus X_2 \oplus \cdots \oplus X_r)T
$$
with $X_j$ commuting with $C_j$, and $X_j$ a block upper triangular Toeplitz matrix
\begin{equation}\label{eq:XblockToeplitzform}
X_j=\begin{bmatrix}  X_{1,j} & X_{2,j} & \cdots &X_{d_j,j}\\ 0 & \ddots & \ddots & \vdots \\ \vdots &\ddots & \ddots & X_{2,j} \\ 0 & \cdots & 0 & X_{1,j}\end{bmatrix}.
\end{equation}
\end{proposition}

\noindent
{\bf Proof.}
Since $A$ is nonderogatory we have that the polynomials $g_j(\lambda)$ are coprime.
Using Theorem \ref{thm:companion-Jordanform} we have
$$A=T^{-1} (C_1\oplus C_2\oplus\cdots \oplus C_r)T,$$ with each $C_j$ given by \eqref{eq:companion-Jordanform}.

The proof proceeds in several steps.

{\it Step 1.}
We discuss how a solution $X$ to $p(X)=A$ can be expressed in terms of the companion-Jordan form of $A$. Note that $A$ must commute with $X$. Write $X=T^{-1}\begin{bmatrix} X_{i,j}\end{bmatrix}_{i,j=1}^r T$. The fact that $X$ and $A$ commute gives $C_iX_{i,j}=X_{i,j}C_j$. So $X_{i,j}$ is a solution of $C_iX_{i,j}-X_{i,j}C_j=0$. Now since $A$ is nonderogatory, $C_i$ and $C_j$ have no common eigenvalues for $i\not= j$. Hence, for $i\not= j$ we have $X_{i,j}=0$. For convenience, set $X_i=X_{i,i}$.
Thus any solution $X$ to $p(X)=A$ must be of the form 
$$
X=T^{-1} (X_1\oplus X_2 \oplus \cdots \oplus X_r)T
$$
with $X_j$ such that $X_jC_j=C_jX_j$, and $p(X_j)=C_j$.

{\it Step 2.}
If $g(\lambda)$ is a monic irreducible polynomial of degree $k$ over $\mathbb{H}$, then $C_g$ is a simple matrix. Then there is a complex $k\times k$ invertible matrix $S$  such that $C_gS=SD$ where $D$ is a diagonal matrix with the eigenvalues of $C_g$ as entries. Let $S$ be the $k\times k$ Vandermonde matrix corresponding to the eigenvalues $\lambda_1,\ldots,\lambda_k$ (the roots of $g(\lambda)$):
\begin{equation*}
S=\begin{bmatrix}
1  & 1& \cdots &  1 \\
\lambda_1&\lambda_2 &\cdots&\lambda_k \\
\lambda_1^2& \lambda_2^2&\cdots&\lambda_k^2\\
\vdots&\vdots&&\vdots\\
\lambda_1^{k-1}&\lambda_2^{k-1}&\cdots&\lambda_k^{k-1}
\end{bmatrix}.
\end{equation*}
Then 
\begin{equation}\label{eq:similarityCD}
\begin{small}
\begin{bmatrix}
S^{-1}&&&\\
&\ddots&&\\
&&\ddots&\\
&&&S^{-1}
\end{bmatrix}\begin{bmatrix}
C_g & I &&\\
&\ddots&\ddots&\\
&&\ddots&I\\
&&&C_g
\end{bmatrix}\begin{bmatrix}
S&&&\\
&\ddots&&\\
&&\ddots&\\
&&&S
\end{bmatrix}=\begin{bmatrix}
D & I &&\\
&\ddots&\ddots&\\
&&\ddots&I\\
&&&D
\end{bmatrix}. 
\end{small}
\end{equation}
Let $Y$ be the matrix similar to $X$ as follows:
\begin{equation}\label{eq:similarityXY}
\begin{bmatrix}
S^{-1}&&\\
&\ddots&\\
&&S^{-1}
\end{bmatrix}X\begin{bmatrix}
S&&\\
&\ddots&\\
&&S
\end{bmatrix}=Y.
\end{equation}
Therefore the solution $X$ of $p(X)=A$ is of the desired block-Toeplitz form as in \eqref{eq:XblockToeplitzform} if and only if $Y$ has block-Toeplitz form. Since $A$ is a polynomial in $X$, it is easy to see that $A$ and $X$ commute. Then from \eqref{eq:similarityCD} and \eqref{eq:similarityXY}, we also know that $Y$ commutes with
\begin{equation}\label{eq:D_ellMatrix}
D_d=\begin{bmatrix}
D&I&0&\cdots&0\\
0&D&I&\ddots&\vdots\\
\vdots&\ddots&\ddots&\ddots&0\\
\vdots&&\ddots&D&I\\
0&\cdots&\cdots&0&D
\end{bmatrix},
\end{equation}
where $d$ is the number of block matrices $D$ on the diagonal. 

Thus it suffices to prove that since $Y$ commutes with the matrix $D_d$,  $Y$  must be a block upper triangular Toeplitz matrix:
\begin{equation*}
Y=\begin{bmatrix}
 Y_1 & Y_2 & \cdots &Y_d\\ 
 0 & \ddots & \ddots & \vdots \\ 
 \vdots &\ddots & \ddots & Y_2 \\ 
 0 & \cdots & 0 & Y_1
 \end{bmatrix}.
\end{equation*}
Moreover, we prove that each $Y_i$ is a diagonal matrix. In fact, the eigenvalues, and consequently the diagonal entries of $D$, are all distinct since the characteristic polynomial of $D$ is the irreducible polynomial $g(\lambda)$. Therefore $Y_iD=DY_i$ implies $Y_i$ is diagonal.

{\it Step 3.1.} Next we proceed by induction on $n$.
Firstly, let $n=2$, that is, $Y$ and $D$ are $2\times 2$ block matrices. We show that if $Y$ commutes with 
\begin{equation*}
D_2=\begin{bmatrix}
D&I\\0&D
\end{bmatrix},
\end{equation*}
then \begin{equation}\label{eq:Y2x2correctform}
Y=\begin{bmatrix}
Y_1&Y_2\\0&Y_1
\end{bmatrix}
\end{equation}
where $Y_1$ and $Y_2$ are diagonal. Let $YD_2=D_2Y$ and write
\begin{equation*}
Y=\begin{bmatrix}
Y_{11}&Y_{12}\\
Y_{21}&Y_{22}
\end{bmatrix}.
\end{equation*}
By comparing the $(2,1)$ blocks in $YD_2=D_2Y$ we have $Y_{21}D=DY_{21}$. Since $Y_{21}$ commutes with $D$, we know $Y_{21}$ is also a diagonal matrix. Next, compare the $(1,1)$ blocks, then we have $Y_{11}D=DY_{11}+Y_{21}$. Consider the $(i,j)$ entries of the matrices in this equation:
\begin{equation*}
(Y_{11})_{ij}\lambda_j=\lambda_i(Y_{11})_{ij}+(Y_{21})_{ij}.
\end{equation*}
If $i\neq j$, then $Y_{21}$ is diagonal, and so $(Y_{11})_{ij}\lambda_j=\lambda_i(Y_{11})_{ij}$. Also, since $\lambda_i\neq\lambda_j$, then $(Y_{11})_{ij}=0$ which means that $Y_{11}$ is also a diagonal matrix. If $i=j$, then
\begin{equation*}
(Y_{11})_{jj}\lambda_j=\lambda_j(Y_{11})_{jj}+(Y_{21})_{jj},
\end{equation*} 
which implies that the diagonal entries of $Y_{21}$ are all zero, thus $Y_{21}=0$.

Comparing the $(2,2)$ blocks in $YD_2=D_2Y$, gives $Y_{21}+Y_{22}D=DY_{22}$. Since $Y_{21}$ is zero, the matrix $Y_{22}$ commutes with $D$, thus $Y_{22}$ is also diagonal. Lastly, we compare the $(1,2)$ blocks and find $Y_{11}+Y_{12}D=DY_{12}+Y_{22}$. Writing $Y_{12}D=DY_{12}+(Y_{22}-Y_{11})$ and following the same argument as with the $(1,1)$ blocks, we obtain that $Y_{12}$ is diagonal and $Y_{22}-Y_{11}=0$, that is, $Y_{22}=Y_{11}$. Therefore $Y$ is of the form \eqref{eq:Y2x2correctform} and all the blocks are diagonal matrices.

{\it Step 3.2.}
Secondly, suppose that for a specific $\nu \geq 2$ any $\nu\times \nu$ block matrix $Y$ commuting with $D_\nu$ has the form
\begin{equation}\label{eq:Y_ellcorrectform}
Y=\begin{bmatrix}
 Y_1 & Y_2 & \cdots &Y_\nu\\ 
 0 & \ddots & \ddots & \vdots \\ 
 \vdots &\ddots & \ddots & Y_2 \\ 
 0 & \cdots & 0 & Y_1
 \end{bmatrix},
\end{equation}
where each $Y_i$ is diagonal.

We use this induction assumption and prove that the above statements hold for $\nu +1$. Let $\widetilde{Y}$ commute with
\begin{equation*}
D_{\nu+1}=\left[\begin{array}{cccc|c}
D&I&&&0\\
&\ddots&\ddots&&\vdots\\
&&\ddots&I&0\\
&&&D&I\\\hline
0&\cdots&\cdots&0&D
\end{array}\right],
\end{equation*}
and suppose the blocks in $\widetilde{Y}$ are given as follows
\begin{equation*}
\widetilde{Y}=\left[\begin{array}{c|c}
Y_{\nu\times\nu}&Z^{(2)}\\ \hline
Z^{(1)}&Z^{(3)}
\end{array}\right],
\end{equation*}
where $Y_{\nu\times\nu}$ is an $\nu\times\nu$ block matrix,
\begin{equation*}
Z^{(1)}=\begin{bmatrix}
Z_1^{(1)}&Z_2^{(1)}&\cdots&Z_\nu^{(1)}
\end{bmatrix}
\end{equation*}
and
\begin{equation*}
Z^{(2)}=\begin{bmatrix}
Z_1^{(2)}\\Z_2^{(2)}\\ \vdots \\ Z_\nu^{(2)}
\end{bmatrix}.
\end{equation*}
All we want to prove, is the following  
\begin{equation*}
Z^{(1)}=0,\;Z^{(2)}=\begin{bmatrix}
Y_{\nu+1}\\Y_\nu \\ \vdots \\ Y_2
\end{bmatrix},\; Z^{(3)}=Y_1,\;  \mbox{\ and\ }Y_{\nu\times\nu}=\begin{bmatrix} 
Y_1 & Y_2 & \cdots &Y_\nu\\
&\ddots&\ddots&\vdots\\
&&\ddots & Y_2\\
&&&Y_1
\end{bmatrix},
\end{equation*}
where $Y_{\nu+1}$ is a diagonal matrix.

As before, we will compare certain block entries in the matrices $\widetilde{Y}D_{\nu+1}$ and $D_{\nu+1}\widetilde{Y}$. We start with block entry $(\nu+1,k)$ for $k=1,\ldots,\nu$. From the $(\nu+1,1)$ blocks we get that $Z_1^{(1)}D=DZ_1^{(1)}$, implying that $Z_1^{(1)}$ is diagonal. Considering the $(\nu+1,2)$ blocks we have 
\begin{equation}\label{eq:ell+1,2}
Z_1^{(1)}+Z_2^{(1)}D=DZ_2^{(1)}
\end{equation}
and comparing the $(\nu+1,k)$ blocks for $k=1,\ldots,\nu$, we have
\begin{equation}\label{eq:ell+1,k}
Z_{k-1}^{(1)}+Z_{k}^{(1)}D=DZ_k^{(1)}.
\end{equation}
Taking into account that $Z_1^{(1)}$ is diagonal and $\lambda_i\neq \lambda_j$ if $i\neq j$, then by considering off-diagonal entries in the blocks in \eqref{eq:ell+1,2}, we obtain that $Z_2^{(1)}$ is diagonal and by considering diagonal entries in the blocks, we find that $Z_1^{(1)}=0$. Doing the same for \eqref{eq:ell+1,k} for each $k$, yields $Z_{k-1}^{(1)}=0$ for $k=1,\ldots,\nu$ and $Z_\nu^{(1)}$ is diagonal.

Next, we consider the blocks $(\nu+1,\nu+1)$ in $\widetilde{Y}D_{\nu+1}=D_{\nu+1}\widetilde{Y}$ and obtain that $Z_\nu^{(1)}=0$ and $Z^{(3)}$ is diagonal.

From the induction assumption we know that $Y_{\nu\times\nu}$ is of the form \eqref{eq:Y_ellcorrectform} since if $\widetilde{Y}$ commutes with $D_{\nu+1}$ and $Z^{(1)}=0$, then $Y_{\nu\times\nu}$ commutes with $D_\nu$.  Hence, we can write
\begin{equation*}
\widetilde{Y}=\left[\begin{array}{cccc|c}
Y_1 & Y_2 & \cdots &Y_\nu &Z_1^{(2)}\\
&\ddots&\ddots&\vdots&Z_2^{(2)}\\
&&\ddots & Y_2&\vdots\\
&&&Y_1&Z_\nu^{(2)}\\\hline
0&\cdots&0&0&Z^{(3)}
\end{array}\right].
\end{equation*}

Finally, compare the $(k,\nu+1)$ blocks in $\widetilde{Y}D_{\nu+1}=D_{\nu+1}\widetilde{Y}$ for\\ $k=\nu,\nu-1,\ldots,1$. For $k=\nu$, we have 
$$Y_1+Z_\nu^{(2)}D=DZ_\nu^{(2)}+Z^{(3)}.$$
Then consider the $(i,j)$ entries of the matrices in this equation and remember that $Y_1$ and $Z^{(3)}$ are diagonal:
\begin{eqnarray*}
\text{if }i\neq j&:&\left(Z_\nu^{(2)}\right)_{ij}\lambda_j=\lambda_i\left(Z_\nu^{(2)}\right)_{ij}\\
\text{if }i=j&:& (Y_1)_{jj}+\left(Z_\nu^{(2)}\right)_{jj}\lambda_j=\lambda_j\left(Z_\nu^{(2)}\right)_{jj}+\left(Z^{(3)}\right)_{jj}.
\end{eqnarray*}
From this we obtain that $Z_\nu^{(2)}$ is diagonal and $Z^{(3)}=Y_1$. Continuing in the same way for $k=\nu-1$ until $k=1$, we use
$$Y_{\nu-k+1}+Z_k^{(2)}D=DZ_k^{(2)}+Z_{k+1}^{(2)}$$
and when comparing the off-diagonal entries in the matrices, one finds that $Z_k^{(2)}$ is diagonal and equating the diagonal entries yields that $Z_{k+1}^{(2)}=Y_{\nu-k+1}$. This completes the proof. \hfill $\Box$

\section{Drazin's results}\label{sec:Drazin}

In this section we present a summary of several results of \cite{Drazin} and discuss the construction given in \cite{Drazin}.

\begin{proposition}\label{prop:Drazin}
Suppose $A$ is nonderogatory with entries in $\mathbb{H}$, and $p(\lambda)$ is a polynomial with coefficients in $\mathbb{H}$. Then the following statements hold:
\begin{itemize}
\item[{\rm  (i)}]
If $p(X)=A$ for some $X$ with entries in $\mathbb{H}$, then $p(\mu)=\lambda_0$ has a solution $\mu\in\mathbb{H}(\lambda_0)$ for every eigenvalue $\lambda_0$ of $A$.
\item[{\rm (ii)}]
Suppose $A$ is simple, and assume that $p(\mu)=\lambda_0$ has a solution $\mu\in\mathbb{H}(\lambda_0)$ for every eigenvalue $\lambda_0$ of $A$. Then there is a solution $X$ with entries in $\mathbb{H}$ of $p(X)=A$.
\item[{\rm (iii)}]
In addition, such solutions $X$ can be constructed explicitly.
\end{itemize}
\end{proposition}

The construction mentioned in the last part of the proposition is as follows.
Let $g(\lambda)$ be an irreducible polynomial over a field $\mathbb{H}$, and suppose that $g(\lambda)$ is the characteristic polynomial of an $n\times n$ matrix $A$. Then $A$ is a simple matrix, as each irreducible polynomial cannot have multiple roots (see for example Pinter \cite{pinter} Chapter 31, Theorem 1). Write $A=U C_{g} U^{-1}$, where $C_{g}$, as usual, is the companion matrix of $g(\lambda)$. Observe, if $Au_j=\lambda_j u_j$, with $u_j\not=0$, then $C_{g} U^{-1}u_j=\lambda_j U^{-1}u_j$, so $U^{-1}u_j$ is an eigenvector of $C_{g}$ corresponding to the eigenvalue $\lambda_j$. By scaling we may then assume that 
$U^{-1} u_j =v_n(\lambda_j)$, i.e., $u_j=Uv_n(\lambda_j)$. We note in passing that this observation can be used to provide an alternative proof of Theorem 3.1 in Drazin's paper~\cite{Drazin}.

Let $\mathbb{H}, \mathbb{K}$ be two fields such that $\mathbb{H}$ is a subfield of $\mathbb{K}$, with $\mathbb{K}$ algebraically closed. 
We now state the main results of Theorem 4.1 in Drazin's paper \cite{Drazin}, as it applies to the situation of the previous paragraph, i.e., $A$ is a simple $n\times n$ matrix over $\mathbb{H}$ with an irreducible characteristic polynomial $g(\lambda)$. Let $p(\lambda)\in\mathbb{H}[\lambda]$. For each eigenvalue $\lambda_j$ of $A$, let $m$ denote the number of different solutions $\mu\in \mathbb{H}(\lambda_j)$ of the equation $p(\mu)=\lambda_j$. Note that $m$ is independent of the selected roots $\lambda_j$ of the characteristic polynomial $g(\lambda)$. Then the equation $p(X)=A$ has exactly $m$ solutions $X$ over $\mathbb{H}$, and these can be constructed explicitly as follows. Let $\mu_j$ be a solution of $p(\mu_j)=\lambda_j$.
Note that $\mu_j u_j$ is a vector with coordinates in $\mathbb{H}(\lambda_j)$, and so we can write
$$
\mu_j u_j=c_{j,0}+c_{j,1}\lambda_j +\cdots +c_{j,n-1}\lambda_j^{n-1}
$$
for some vectors $c_{j,i}$ with coordinates in $\mathbb{H}$. 
Define the matrix 
$$M=\begin{bmatrix} c_{j,0} & c_{j,1} &\cdots & c_{j,n-1}\end{bmatrix}.$$ 
Then, by Proposition 2.4 in \cite{Drazin}, if $p(X)=A$ for some matrix $X$ with entries in $\mathbb{H}$ we must have $XU=M$, in other words $X=MU^{-1}$, and conversely, by Theorem 4.1 in \cite{Drazin}, if $X$ is constructed this way, it will be a solution to $p(X)=A$.

\section{Solving $p(X)=A$ for the case of nonderogatory $A$}\label{sec:solving}

As explained in the previous sections we will assume from the start that the matrix $A$ is given in the following $dk\times dk$  companion-Jordan form 
$$
A=\begin{bmatrix} 
C_g & I_k & 0 & \cdots & 0\\
0 & C_g & I_k & \ddots & \vdots \\
\vdots & \ddots & \ddots & \ddots &0\\
\vdots  & \cdots & \ddots & C_g& I_k\\
0 & \cdots & \cdots & 0 & C_g
\end{bmatrix},
$$
with $g(\lambda)=g_0 + g_1\lambda+ \cdots +g_{k-1}\lambda^{k-1}+g_k\lambda^k$ a monic (so $g_k=1$) and irreducible polynomial of degree $k$ over $\mathbb{H}$.
For a monic polynomial $p(\lambda)=p_0+p_1\lambda+\cdots +p_{\ell-1}\lambda^{\ell-1}+\lambda^\ell\in\mathbb{H}[\lambda]$ of degree $\ell$, we shall look for a solution of $p(X)=A$ of the form
$$
X=\begin{bmatrix}  X_1 & X_2 & \cdots &X_d\\ 0 & \ddots & \ddots & \vdots \\ \vdots & & \ddots & X_2 \\ 0 & \cdots & 0 & X_1\end{bmatrix}.
$$
With such a matrix $X$ we have that $p(X)$ is upper triangular Toeplitz, with diagonal block $p(X_1)$. It follows that $X_1$ satisfies $p(X_1)=C_g$. From the previous section, in fact, from \cite{Drazin}, we know what the requirements are for a solution $X_1$ to exist, and an explicit construction is given there as well.

Observe that if $p(X)=A$, then not only $p(X_1)=C_g$, but also 
$$
p\left(\begin{bmatrix} X_1 & X_2 \\ 0 & X_1\end{bmatrix}\right) = 
\begin{bmatrix} C_g & I_k \\ 0 & C_g\end{bmatrix}.
$$
From the results in Section \ref{sec:deltas} we have
$$
p\left(\begin{bmatrix} X_1 & X_2 \\ 0 & X_1\end{bmatrix} \right) =\begin{bmatrix} p(X_1) & \Delta p(X_1,X_1)(X_2)\\ 0 & p(X_1)\end{bmatrix}.
$$
So, $X_2$ may be determined from the linear matrix equation
\begin{equation*}
 \Delta p(X_1,X_1)(X_2) =I_k.
\end{equation*}
To find $X_2$, solve (when possible) the linear matrix equation
\begin{equation}\label{eq:X2}
\Delta p (X_1,X_1)(X_2)=\sum_{m=1}^{\ell} p_m\sum_{j=0}^{m-1} X_1^{m-1-j}X_2X_1^j =I_k.
\end{equation}
Now this equation may or may not be solvable, but if it is solvable over the algebraically closed field $\mathbb{K}$, then it is solvable over $\mathbb{H}$, by the results in Section~\ref{sec:lineq}.

Next, to find $X_3$ we consider the $3d\times 3d$ block in the upper right hand corner, and the equation
$$
p\left(\begin{bmatrix} X_1 & X_2&X_3 \\ 0 & X_1& X_2 \\ 0 & 0 & X_1\end{bmatrix} \right) =\begin{bmatrix} C_g & I_k & 0 \\ 0 & C_g & I_k \\ 0 & 0 & C_g\end{bmatrix}.
$$
Partition the matrices in the equation as follows
$$p\left(\left[\begin{array} {cc|c}
X_1 & X_2 & X_3 \\ 0 & X_1 & X_2\\ \hline0 & 0 & X_1\end{array}\right]
\right) =\left[\begin{array} {cc|c} C_g & I_k & 0  \\ 0 & C_g& I_k \\ \hline 0 & 0 & C_g\end{array}\right].
$$
Equivalently we have to
solve the following linear matrix equation for $X_3$:
{\small
$$
\Delta p\left(\begin{bmatrix} X_1 & X_2\\ 0 & X_1\end{bmatrix}, X_1\right)\left(\begin{bmatrix} X_3 \\ X_2 \end{bmatrix}\right)= 
\sum_{m=1}^{\ell} p_m\sum_{j=0}^{m-1} 
\begin{bmatrix} X_1 & X_2\\ 0 & X_1\end{bmatrix}^{m-1-j}
\begin{bmatrix} X_3 \\ X_2 \end{bmatrix}
X_1^j =\begin{bmatrix} 0 \\ I_k \end{bmatrix}.
$$
}

Continue in this manner to solve consecutively for $X_4, \ldots , X_d$.

Since each step gives a set of linear equations in the entries in the matrix $X_j$ to be solved, it can be shown inductively that if a solution exists it must have entries in $\mathbb{H}$. 

To summarize, the algorithm becomes: 
\\
Step 1. Find $X_1$ (using the result from Section \ref{sec:Drazin}).
\\
Step 2. Solve $X_2$ from
$$
\Delta p (X_1,X_1)(X_2)=\sum_{m=1}^{\ell} p_m\sum_{j=0}^{m-1} X_1^{m-1-j}X_2X_1^j =I_k.
$$
\\
Step 3. For $j=3, 4, \ldots, d$ solve consecutively for $X_j$ the linear matrix equation
$$
\Delta p\left(\begin{bmatrix} 
X_1 & X_2 & \cdots & \cdots  & X_{j-1}\\ 
0 & X_1 & \ddots &   & \vdots \\
\vdots & \ddots &  \ddots & \ddots & \vdots\\
\vdots &  & \ddots & X_1 & X_2 \\
0 & \cdots  & \cdots & 0 & X_1
\end{bmatrix}, X_1\right)
\left( \begin{bmatrix} X_j \\ X_{j-1} \\ \vdots \\ \vdots \\ X_2\end{bmatrix} \right)= 
\begin{bmatrix} 0 \\ \vdots \\ \vdots \\ 0 \\  I_k \end{bmatrix}.
$$

If all these equations have a solution, then $p(X)=A$ has a solution in~$\mathbb{H}^{n\times n}$.

Let us analyze the equations a bit more. Introduce
$$
Z=\begin{bmatrix} X_2 & \cdots & X_{j-1}\end{bmatrix},\  \widetilde{X}_1=\begin{bmatrix}
X_1 & X_2 & \cdots & \cdots  & X_{j-2}\\ 
0 & X_1 & \ddots &   & \vdots \\
\vdots & \ddots &  \ddots & \ddots & \vdots\\
\vdots &  & \ddots & X_1 & X_2 \\
0 & \cdots  & \cdots & 0 & X_1
\end{bmatrix}, \ \widetilde{X}=\begin{bmatrix} X_{j-1}\\ \vdots\\ \vdots \\ \vdots \\ X_2\end{bmatrix}.
$$
Then the equation
\begin{equation}\label{eq:Xj}
\Delta p\left(\begin{bmatrix} 
X_1 & X_2 & \cdots & \cdots  & X_{j-1}\\ 
0 & X_1 & \ddots &   & \vdots \\
\vdots & \ddots &  \ddots & \ddots & \vdots\\
\vdots &  & \ddots & X_1 & X_2 \\
0 & \cdots  & \cdots & 0 & X_1
\end{bmatrix}, X_1\right)
\left( \begin{bmatrix} X_j \\ X_{j-1} \\ \vdots \\ \vdots \\ X_2\end{bmatrix} \right)= 
\begin{bmatrix} 0 \\ \vdots \\ \vdots \\ 0 \\  I_k \end{bmatrix}
\end{equation}
becomes
$$
\sum_{m=1}^\ell p_m \sum_{i=0}^{m-1} \begin{bmatrix} X_1 & Z\\ 0 & \widetilde{X}_1\end{bmatrix}^{m-1-i}
\begin{bmatrix} X_j \\ \widetilde{X}\end{bmatrix} X_1^i=\begin{bmatrix} 0 \\ \vdots \\ 0 \\ I_k\end{bmatrix}.
$$
Now 
$$
\begin{bmatrix} X_1 & Z\\ 0 & \widetilde{X}_1\end{bmatrix}^{m-1-i} =\begin{bmatrix} X_1^{m-1-i} & W_i \\ 0 &  \widetilde{X}_1^{m-1-i}\end{bmatrix}
$$
for some $W_i$ which depends only on $X_1, \ldots , X_{j-1}$.
Thus \eqref{eq:Xj} becomes
{\small
$$
\sum_{m=1}^\ell p_m \sum_{i=0}^{m-1} \begin{bmatrix} 
X_1^{m-1-i} & W_i \\ 0 &  \widetilde{X}_1^{m-1-i}\end{bmatrix}
\begin{bmatrix} X_j \\ \widetilde{X}\end{bmatrix} X_1^i=\sum_{m=1}^\ell p_m \sum_{i=0}^{m-1} 
\begin{bmatrix} X_1^{m-1-i}X_jX_1^i+W_i\widetilde{X}X_1^i \\
\widetilde{X}_1^{m-1-i}\widetilde{X}X_1^i \end{bmatrix}
=
\begin{bmatrix} 0 \\ \vdots \\ 0 \\ I_k\end{bmatrix}.
$$
}
Observe that $\sum_{m=1}^\ell p_m \sum_{i=0}^{m-1} \widetilde{X}_1^{m-1-i}\widetilde{X}X_1^i 
=\Delta p(\widetilde{X}_1, X_1)(\widetilde{X})$, which we know already from step $j-1$ is equal to $\begin{bmatrix} 0 \\ \vdots \\ 0 \\ I_k\end{bmatrix}$. So only the top component gives an extra identity from which we need to solve $X_j$.

Hence equation \eqref{eq:Xj} is equivalent to
$$
0=\sum_{m=1}^\ell p_m \sum_{i=0}^{m-1} X_1^{m-1-i}X_jX_1^i 
+ \sum_{m=1}^\ell p_m \sum_{i=0}^{m-1}W_i\widetilde{X}X_1^i,
$$
or alternatively,
$$
\sum_{m=1}^\ell p_m \sum_{i=0}^{m-1} X_1^{m-1-i}X_jX_1^i  = -\sum_{m=1}^\ell p_m \sum_{i=0}^{m-1}W_i\widetilde{X}X_1^i.
$$

Next, we compute the $W_i$ explicitly. When $m-1-i=0$, so $i=m-1$, we have $W_i=0$, when $i=m-2$ we have $W_i=Z$, when $i=m-3$ we have $W_i=X_1Z+Z\widetilde{X}_1$ and when $i=m-4$ we have $W_i=X_1^2Z+X_1Z\widetilde{X}_1 +Z\widetilde{X}_1^2$. We now see a pattern emerging. It can be shown that we have 
$$
W_i=\Delta_{\lambda^{m-1-i}}(X_1,\widetilde{X}_1)(Z)=\sum_{\nu=0}^{m-i-2}X_1^{m-i-2-\nu}\begin{bmatrix} X_2 & \cdots & X_{j-1}\end{bmatrix}\widetilde{X}_1^\nu.
$$

So, also using $W_{m-1}=0$, \eqref{eq:Xj} becomes,
\begin{align*}
&\sum_{m=1}^\ell p_m \sum_{i=0}^{m-1} X_1^{m-1-i}X_jX_1^i   \\ =&-\sum_{m=1}^\ell p_m \sum_{i=0}^{m-2} \sum_{\nu=0}^{m-i-2}X_1^{m-i-2-\nu}\begin{bmatrix} X_2 & \cdots & X_{j-1}\end{bmatrix}\widetilde{X}_1^\nu
\begin{bmatrix} X_{j-1}\\ \vdots \\ X_2\end{bmatrix} X_1^i,
\end{align*}
with $\widetilde{X}_1 =\begin{bmatrix}
X_1 & X_2 & \cdots & \cdots  & X_{j-2}\\ 
0 & X_1 & \ddots &   & \vdots \\
\vdots & \ddots &  \ddots & \ddots & \vdots\\
\vdots &  & \ddots & X_1 & X_2 \\
0 & \cdots  & \cdots & 0 & X_1
\end{bmatrix}$.

Applying the Kronecker product to rewrite the matrix equation, we now obtain
\begin{align*}
& \left(\sum_{m=1}^\ell p_m \sum_{i=0}^{m-1} (X_1^T)^i \otimes X_1^{m-1-i}\right) \Vec{X_j}
\nonumber \\
=&
-\Vec{\sum_{m=1}^\ell p_m \sum_{i=0}^{m-2} \sum_{\nu=0}^{m-i-2}X_1^{m-i-2-\nu}\begin{bmatrix} X_2 & \cdots & X_{j-1}\end{bmatrix}\widetilde{X}_1^\nu
\begin{bmatrix} X_{j-1}\\ \vdots \\ X_2\end{bmatrix} X_1^i}.
\end{align*}
Note that on the right hand side the $m=1$ term is always an empty sum, so that term vanishes.
Hence the equation finally becomes
\begin{align}\label{eq:XjKron}
& \left(\sum_{m=1}^\ell p_m \sum_{i=0}^{m-1} (X_1^T)^i \otimes X_1^{m-1-i}\right) \Vec{X_j}
\nonumber \\
=&
-\Vec{\sum_{m=2}^\ell p_m \sum_{i=0}^{m-2} \sum_{\nu=0}^{m-i-2}X_1^{m-i-2-\nu}\begin{bmatrix} X_2 & \cdots & X_{j-1}\end{bmatrix}\widetilde{X}_1^\nu
\begin{bmatrix} X_{j-1}\\ \vdots \\ X_2\end{bmatrix} X_1^i}.
\end{align}

{\bf Remarks} (1) The coefficient matrix $\sum_{m=1}^\ell p_m \sum_{i=0}^{m-1} (X_1^T)^i \otimes X_1^{m-1-i}$ is independent of $j$, whereas the right hand side of \eqref{eq:XjKron} depends on $j$. This means that the following holds.
Let $p(\lambda)$ be such that for each root $\lambda_0$ of $g(\lambda)$ (so each eigenvalue of $A$) there is a $\mu\in\mathbb{H}(\lambda_0)$ such that $p(\mu)=\lambda_0$, and let $X_1\in\mathbb{H}^{k\times k}$ be such that $p(X_1)=C_g$. If 
\begin{equation}\label{eq:detcondition}
\det \left(\sum_{m=1}^\ell p_m \sum_{i=0}^{m-1} (X_{1}^T)^i \otimes X_{1}^{m-1-i}\right)\not=0
\end{equation}
then equation \eqref{eq:XjKron} is uniquely solvable for each $j$, and so $X_2, \ldots , X_d$ depend uniquely on $X_1$.\hfill $\Box$

\bigskip

(2)
Equation \eqref{eq:X2} can also be rewritten in terms of the Kronecker product as follows: 
$\left(\sum_{m=1}^\ell p_m \sum_{i=0}^{m-1} (X_1^T)^i \otimes X_1^{m-1-i}\right)\Vec{X_2}=\Vec{I_k}$. The condition for unique solvability is again \eqref{eq:detcondition}.

Note that we have to be a bit careful here: it could still be the case that equations \eqref{eq:X2} and \eqref{eq:XjKron} have a solution (and then in fact infinitely many) when $\det\left(\sum_{m=1}^\ell p_m \sum_{i=0}^{m-1} (X_1^T)^i \otimes X_1^{m-1-i}\right)=0$.
\hfill$\Box$

\bigskip

(3) Note that in the scalar case the condition \eqref{eq:detcondition} becomes 
$$p^\prime(X_{1})\not=0.$$
Further, equation \eqref{eq:X2} becomes $p^\prime(X_{1})X_{2}=1$, and for $k=3$ the equation \eqref{eq:XjKron} becomes (after some rearrangement and grouping of terms) $p^\prime(X_{1})X_{3}=-\tfrac{1}{2}p^{\prime\prime}(X_{1})X_{2}^2$. For $k=4$ and higher, the equations in \eqref{eq:XjKron} become less transparent even in the scalar case. \hfill $\Box$

\bigskip

(4)
Note that the relation goes deeper than just the scalar case. Let us assume \eqref{eq:detcondition} is satisfied. Let $\mu$ be an eigenvalue of $X_{1}$, and let $X_{1}x=\mu x$, and $X_{1}^T y=\mu y$. Then for every $i$ and $m$ we have 
$$
\left((X_{1}^T)^i \otimes X_{1}^{m-1-i}\right) (y\otimes x) =(X_{1}^T)^i y \otimes X_{1}^{m-1-i}x =\mu^i y\otimes \mu^{m-i-1}x=\mu^{m-1} (y\otimes x).
$$
Hence
\begin{align*}
& \left(\sum_{m=1}^\ell p_m \sum_{i=0}^{m-1} (X_{1}^T)^i \otimes X_{1}^{m-1-i}\right) (y\otimes x)= \left(\sum_{m=1}^\ell p_m\sum_{i=0}^{m-1} \mu^{m-1}\right) (y\otimes x)\\  =&
\left(  \sum_{m=1}^\ell mp_m  \mu^{m-1}\right) (y\otimes x) =p^\prime (\mu)  (y\otimes x).
\end{align*}
So $p^\prime(\mu)$ is an eigenvalue of $\sum_{m=1}^\ell p_m \sum_{i=0}^{m-1} (X_{1}^T)^i \otimes X_{1}^{m-1-i}$ with corresponding eigenvector $y\otimes x$.
If condition \eqref{eq:detcondition} holds, then it follows that $p^\prime(\mu)\not=0$. 

We see that a necessary condition for \eqref{eq:detcondition} to hold is that $p^\prime(\mu)\not=0$ for the eigenvalues $\mu$ of $X_{1}$. 

This has also a necessary and sufficient version. Using the result of exercise 4.2.19 in \cite{HornJohn}, we have the following. Introduce the polynomial in two variables:
$$
\widetilde{p}(s,t)= \sum_{m=1}^\ell p_m \sum_{i=0}^{m-1} s^it^{m-1-i}.
$$
(Observe that $\widetilde{p}(s,s)=p^\prime (s)$.)
Then
\begin{align*}
\sigma\left(\sum_{m=1}^\ell p_m \sum_{i=0}^{m-1} (X_{1}^T)^i \otimes X_{1}^{m-1-i}\right) &=\left\{ \widetilde{p}(\mu_1,\mu_2) \mid \mu_1, \mu_2\in\sigma(X_{1})\right\}\\
&=\left\{  \sum_{m=1}^\ell p_m \sum_{i=0}^{m-1} \mu_1^i \mu_2^{m-1-i} \mid \mu_1, \mu_2\in\sigma(X_{1})\right\}.
\end{align*}
In particular, condition \eqref{eq:detcondition} is equivalent to 
\begin{equation}\label{eq:eigenvaluecondition}
\sum_{m=1}^\ell p_m \sum_{i=0}^{m-1} \mu_1^i \mu_2^{m-1-i} \not=0
\end{equation}
for each pair $\mu_1, \mu_2$ of eigenvalues of $X_{1}$. \hfill$\Box$

\bigskip

To summarize what we have proved, we state the following theorem, which we consider as the main result of the paper. It is the more precise version of Theorem \ref{thm:main_intro} which was promised in the introduction.
Part (i) is already formulated in \cite{Drazin}, part (ii) follows immediately from \cite{Drazin} (see Proposition \ref{prop:Drazin}). Parts (iii) and (v) are our main contribution. Part (iv) follows from our construction, combined with the number of solutions given in \cite{Drazin} for the simple case.

\begin{theorem}\label{thm:main}
Consider an $n\times n$ nonderogatory matrix $A$ with entries in a field $\mathbb{H}$, with $\mathbb{Q}\subset \mathbb{H}\subset\mathbb{C}$, and let $p_A(\lambda)=g_1(\lambda)^{d_1}\cdot g_2(\lambda)^{d_2}\cdots g_r(\lambda)^{d_r}$ be the factorization of the characteristic polynomial of $A$ with $g_j(\lambda)$'s pairwise coprime monic and irreducible polynomials.
Its companion-Jordan form is given by
\begin{equation*}
A=T^{-1} (C_1\oplus C_2\oplus \cdots \oplus C_r)T
\end{equation*}
where each $C_j$ is a $d_jk_j\times d_jk_j$ block upper triangular matrix of the form
\begin{equation*}
C_j=\begin{bmatrix} 
C_{g_j} & I_{k_j} & 0 & \cdots & 0\\
0 & C_{g_j} & I_{k_j} & \ddots & \vdots \\
\vdots & \ddots & \ddots & \ddots &0\\
\vdots  &  & \ddots & C_{g_j}& I_{k_j}\\
0 & \cdots & \cdots & 0 & C_{g_j}
\end{bmatrix},
\end{equation*}
and where $T\in\mathbb{H}^{n\times n}$ is invertible.

Let $p(\lambda)=p_0+p_1\lambda+\cdots +p_\ell \lambda^\ell \in\mathbb{H}[\lambda]$. 
Then the following hold.
\begin{itemize}
\item[{\rm (i)}]
If $p(X)=A$ has a solution $X$ in $\mathbb{H}^{n\times n}$, then for each root $\lambda_j$ of some $g_k(\lambda)$ there is a solution $\mu\in \mathbb{H}(\lambda_j)$ with $p(\mu)=\lambda_j$. 
\item[{\rm (ii)}] 
Conversely, assume that $p(\mu)=\lambda$ has a solution in $\mathbb{H}(\lambda)$ for each eigenvalue $\lambda$ of $A$. Let $\lambda_j$ be a root of $g_j(\lambda)$, and let $m_j$ denote the number of such solutions of $p(\mu)=\lambda_j$. 
For each $j=1, \ldots, r$ there are $m_j$ solutions $X_{1,j}$  to $p(X_{1,j})=C_{g_j}$. 
\item[{\rm (iii)}]
If for each $j$ for which $d_j>1$ there is at least one solution $X_{1,j}$ for which
\eqref{eq:detcondition} or equivalently \eqref{eq:eigenvaluecondition} is satisfied,
then there is a solution $X$ of $p(X)=A$ with entries in $\mathbb{H}$ and such a solution can be constructed explicitly via the solution of a set of linear equations, as follows. 
The solution is of the form 
$$
X=T^{-1} (Y_1\oplus Y_2 \oplus \cdots \oplus Y_r)T,
$$
where each $Y_j$ is a block upper triangular Toeplitz matrix
$$
Y_j=\begin{bmatrix}
X_{1,j} & X_{2,j} & \cdots & \cdots  & X_{d_j,j}\\ 
0 & X_{1,j} & \ddots &   & \vdots \\
\vdots & \ddots &  \ddots & \ddots & \vdots\\
\vdots &  & \ddots & X_{1,j} & X_{2,j} \\
0 & \cdots  & \cdots & 0 & X_{1,j}
\end{bmatrix}.
$$
Here $X_{1,j}$ is a solution of $p(X_{1,j})=C_{g_j}$, for which \eqref{eq:detcondition} holds. The matrix
$X_{2,j}$ is the unique solution to
\eqref{eq:X2}, 
while for $k=3, \ldots , d_j$ the matrix 
$X_{k,j}$ is the unique solution to
\eqref{eq:XjKron}.
\item[{\rm (iv)}]
If \eqref{eq:detcondition} or equivalently \eqref{eq:eigenvaluecondition} is satisfied
for each $j$ for which $d_j>1$ and for each possible solution $X_{1,j}$, then there are precisely $m_1\cdot m_2 \cdots m_r$ solutions.
\item[{\rm (v)}] 
If  $d_j>1$ and \eqref{eq:detcondition} is not satisfied for some solution $X_{1,j}$ of $p(X_{1,j})=C_{g_j}$, but equations \eqref{eq:X2} and \eqref{eq:XjKron} are solvable with this choice of $X_{1,j}$, then a solution $X\in\mathbb{H}^{n\times n}$ can still be constructed as indicated in part {\rm (iii)} above. In fact, there will be infinitely many solutions in that case.
\end{itemize}
\end{theorem}

\bigskip
We conclude the section with several remarks.

\bigskip

{\bf Remark on the number of solutions.} Part (iv) of the theorem above deals with a situation where there are finitely many solutions. This is not the only situation in which there may be finitely many solutions. After all, it may happen that for one of the possible solutions $X_{1,j}$ originating from part (ii) the condition \eqref{eq:detcondition} is satisfied, while for another the condition is not satisfied, while at the same time equations \eqref{eq:X2} and \eqref{eq:XjKron} are not solvable for the latter solution $X_{1,j}$. In that case there may still be finitely many solutions, and the number mentioned in part (iv) is only an upper bound. 

\bigskip

{\bf Remark on complex solutions.} Note that taking $\mathbb{H}=\mathbb{C}$ we obtain a condition for the existence of a complex solution $X$ of $p(X)=A$ (compare \cite{EU, roth, Spiegel}). However, we can also apply the theorem in a slightly different way. With the notation as in the statement of the theorem, let $\mathbb{H}$ be the smallest field extension of $\mathbb{Q}$ for which the following holds: for each root $\lambda_j$ of some $g_k(\lambda)$ there is a solution $\mu\in \mathbb{H}(\lambda_j)$ with $p(\mu)=\lambda_j$. Suppose in addition that there is a complex solution $X$ of $p(X)=A$. Then there is a solution $X$ of $p(X)=A$ over $\mathbb{H}$.

\bigskip

{\bf Remark on the derogatory case.} In case $A$ is derogatory, we can still use Theorem \ref{thm:companion-Jordanform} to write $A=T^{-1} (C_1 \oplus \cdots \oplus C_r)T$, where each $C_j$ is of the form \eqref{eq:companion-Jordanform}. If for every eigenvalue $\lambda_0$ of $A$ there is a solution $\mu\in\mathbb{H}(\lambda_0)$ of $p(\mu)=\lambda_0$, then one may construct for each  $j$ a solution to $p(X_{1,j})=C_{g_j}$ as in Section \ref{sec:Drazin}. If for each $j$ the set of equations \eqref{eq:X2} and \eqref{eq:XjKron} are solvable for $X_{2,j}, \ldots , X_{d_j}$ then a solution in $\mathbb{H}^{n\times n}$ of $p(X)=A$ exists, and is constructed in the same way as in the nonderogatory case.

However, in the derogatory case a solution may exist that is not of this form. A straightforward example of that is to take $A=I_2$ and $p(\lambda)=\lambda^2$, $\mathbb{H}=\mathbb{Q}$. Then $X=\begin{bmatrix} \tfrac{1}{2} & 1 \\[1mm] \tfrac{3}{4} & -\tfrac{1}{2}\end{bmatrix}$ is a solution. A more fundamental obstruction will be illustrated in the next section: solutions may exist that cannot be constructed per block, but need combinations of several blocks.

\section{Special case $X^\ell=A$}\label{sec:cases}

First consider the case $\ell=2$. 
In this particular case we can be much more explicit, as characterizing the solvability of $\mu^2=\lambda_j$ in the field $\mathbb{H}(\lambda_j)$ is much more straightforward than solving a general polynomial equation. Once that is dealt with, the next steps in the solution of $X^2=A$ are also more straightforward. Indeed, for $p(\lambda)=\lambda^2$ we have that $\Delta p(X_1,X_1)(X_2)=I_k$ becomes the equation
$$
X_1X_2+X_2X_1=I_k.
$$
Hence we can be a bit more specific concerning uniqueness of solvability: if $X_1$ and $-X_1$ do not have common eigenvalues, then there is a unique solution $X_2$ over $\mathbb{C}$, see e.g., Horn and Johnson \cite{HornJohn}, Theorem 4.4.6, p.\ 270, or Lancaster and Tismenetsky \cite{LT}, Theorem 12.3.2, p.\ 414.
Moreover, explicit conditions for existence of a square root of a given matrix $A$ are well-known in terms of eigenvalues and Jordan chains of $A$.

This case also provides us with an example where the second step in the algorithm fails, even when $A$ is nonderogatory. Consider $A=J_2(0)$. In this case $X_1=0$, but the equation $X_1X_2+X_2X_1=1$ obviously has no solution. 

For the derogatory case we have from this also an example where a solution exists, but cannot be obtained from Theorem \ref{thm:main}. Take $A=J_2(0)\oplus (0)$. Then $A=X^2$ for some $X$ which is similar to $J_3(0)$, but there is no way that this solution can be constructed by the algorithm of the theorem.

For the case where $\ell >2$ the condition \eqref{eq:eigenvaluecondition} becomes the following: for every pair of eigenvalues $\mu_1, \mu_2$ of $X_1$ we have
$\sum_{i=0}^{\ell-1} \mu_1^i\mu_2^{\ell-1-i} \not= 0$.

\section{Examples}\label{sec:examples}

\noindent
{\bf Example 1.}
Let us consider the companion matrix $$A=\begin{bmatrix} 0 & 1 & 0 & 0 & 0 & 0  \\ 0 & 0 & 1 & 0& 0 & 0  \\ 0 & 0 & 0 & 1 & 0 & 0 \\ 0 & 0 & 0 & 0 & 1 & 0 \\ 0 & 0 & 0 & 0 & 0 & 1 \\ -64 & 0 & -48 & 0 & -12 & 0 \end{bmatrix},$$ which is nonderogatory and has characteristic polynomial $p_A(\lambda)=(\lambda^2+4)^3$. Our interest is in square roots of $A$, that is, in solving $p(\lambda)=\lambda^2$, particularly in rational square roots.

\subsubsection*{The companion-Jordan form.}
First we bring $A$ to companion-Jordan form. So we have $g(\lambda)=\lambda^2+4$, with $d=3$ and $k=2$.

Take an eigenvector corresponding to the eigenvalue $2i$: 
$$
v_6(2i)=\begin{bmatrix} 1 \\ 2i \\ -4 \\ -8i \\ 16 \\ 32 i \end{bmatrix}=\begin{bmatrix} 1 & 0 \\ 0 & 1 \\ -4 & 0 \\ 0 & -4 \\ 16 & 0 \\ 0 & 16 \end{bmatrix} \begin{bmatrix} 1 \\ 2i \end{bmatrix} =Wv_2(2i).
$$
Then we construct $T=\begin{bmatrix} W & S_1 W & S_2W\end{bmatrix}$, which after some computation equals
$$
T=\begin{bmatrix}1 & 0 & 0 & 0 & 0 & 0  \\ 0 & 1 & 1 & 0 & 0 & 0  \\ -4 & 0 & 0 & 2& 1 & 0  \\ 0 & -4 & -12 & 0 & 0 & 3 \\ 16 & 0 & 0 & -16 & -24 & 0 \\ 0 & 16 & 80 & 0 & 0 & -40\end{bmatrix},
$$
and for this $T$ we have, 
$$
T^{-1}AT=B=\begin{bmatrix} 0 & 1 & 1 & 0 &0 & 0 \\ -4 & 0 & 0 & 1  &0 & 0\\ 0 & 0 & 0 & 1  &1 & 0\\ 0 & 0 & -4 & 0 &0 & 1 \\ 0 & 0 & 0 & 0 & 0 & 1 \\ 0 & 0 & 0 & 0 & -4 & 0\end{bmatrix}
=
\begin{bmatrix} A_1 & I_2 & 0 \\ 0 & A_1 & I_2 \\ 0 & 0 & A_1\end{bmatrix},
$$
which is in companion-Jordan form with $A_1=\begin{bmatrix} 0 & 1 \\ -4 & 0 \end{bmatrix}$.

\subsubsection*{Square roots of the block upper left corner.}
In line with what we did before, consider $A_1=\begin{bmatrix} 0 & 1 \\ -4 & 0 \end{bmatrix}$. We follow the approach by Drazin \cite{Drazin} for
$\mathbb{H}=\mathbb{Q}$
and $\mathbb{K}=\mathbb{C}$, and consider the existence of a rational square root of $A_1$. The characteristic polynomial $p_{A_1}(\lambda)=\lambda^2+4$ has two complex roots $\pm 2i$. Consider $p(\lambda)=\lambda^2$, so we are looking for square roots $X_1$ of $A_1$: $X_1^2=A_1$. Then the solutions to $p(\lambda)=2i$ are $\pm(1+i)$ and the solutions to $p(\lambda)=-2i$ are $\pm (-1+i)$. 

In Theorem 4.1 of Drazin's paper, choose $\lambda_j=2i$. Then $m_j=2$ as there are two different solutions to $p(\mu_j)=\lambda_j$,  namely $\pm (1+i)$, which both are in $\mathbb{Q}(2i)=\mathbb{Q}(i)$. The important point is that the solutions to $p(\lambda)=2i$ are in~$\mathbb{Q}(i)$.

Anyway, let us see how the construction in Theorem 4.1 of Drazin's paper works out. The eigenvector $w$ of $A_1$ corresponding to $2i$ is the vector $\begin{bmatrix} 1\\ 2i \end{bmatrix}$ 
Now $w=w_0+2iw_1$ with
$w_0=\begin{bmatrix} 1 \\ 0\end{bmatrix}$ and $w_1=\begin{bmatrix} 0 \\ 1 \end{bmatrix}$. From equation (9) in Drazin we then have, if we take $\mu=1+i$ that
$
X_1w=(1+i)w
$,
resulting in 
$
X_1\begin{bmatrix} w_0 & w_1\end{bmatrix}= M
$
with $M$ determined by
$$
(1+i)\begin{bmatrix} 1\\ 2i \end{bmatrix} =\begin{bmatrix} 1+i\\ -2 +2i\end{bmatrix}
=\begin{bmatrix} 1 & \tfrac{1}{2}\\ -2 & 1 \end{bmatrix}\begin{bmatrix} 1\\ 2i \end{bmatrix}.
$$
Since in this case $\begin{bmatrix}w_0 & w_1\end{bmatrix}$ happens to be the identity, we get $X_1=\begin{bmatrix} 1 &  \tfrac{1}{2}\\ -2 & 1 \end{bmatrix}$. One easily checks that indeed, this is a solution to $X_1^2=A_1$.

If we had chosen $\mu=-1 -i$ we would have obtained that the corresponding solution $X_1$ must satisfy
$
X_1w=(-1-i)w,
$
yielding 
$
X_1\begin{bmatrix}w_0 & w_1\end{bmatrix}= M
$
with $M$ obtained from
$$
(-1-i)\begin{bmatrix} 1\\ 2i \end{bmatrix} =\begin{bmatrix} -1-i\\ 2 -2i\end{bmatrix}
=\begin{bmatrix} -1 & -\tfrac{1}{2}\\ 2 & -1 \end{bmatrix}\begin{bmatrix} 1\\ 2i \end{bmatrix},
$$
giving the solution $X_1=\begin{bmatrix} -1 & -\tfrac{1}{2}\\ 2 & -1 \end{bmatrix}$.

There are two other complex square roots of $A_1$, namely, $\pm \begin{bmatrix} i & -\tfrac{1}{2}i \\ 2i & i\end{bmatrix}$. Note that these have eigenvalues $1+i$ and $-1+i$, respectively $1-i$ and $-1-i$. These solutions do not come about in the way that Drazin explains, unless we work over $\mathbb{H}=\mathbb{Q}(i)$, as in that field we can factorize $p_A(\lambda)$ as $(\lambda-2i)(\lambda+2i)$. Applying Drazin's Theorem 4.1 to this situation, one finds that there are four solutions (which are indeed all solutions, even over $\mathbb{C}$).

\subsubsection*{Rational square roots of $A$}
Let 
$$
B_1=\begin{bmatrix} A_1 & I_2 \\ 0 & A_1\end{bmatrix}.
$$
Can we now find solutions $Y$ of $Y^2=B_1$ over $\mathbb{Q}$? Let us try a solution of the form
$$
Y=\begin{bmatrix} X_1 & X_{2} \\ 0 & X_1\end{bmatrix},
$$
with $X_1=\begin{bmatrix} 1 & \tfrac{1}{2} \\ -2 & 1 \end{bmatrix}$.
Then we have that $Y^2=B_1$ if and only if 
$$
X_1X_{2}+X_{2}X_1=I_2.
$$
Now since $X_1$ has eigenvalues $1\pm i$ we have $\sigma(X_1)\cap \sigma(-X_1)=\emptyset$, and hence this matrix equation is uniquely solvable. But this matrix equation can also be written as a set of four linear equations in the four entries of $X_{2}$, and then the unique solution can be found using Cramer's rule. Since $X_1$ has entries in $\mathbb{Q}$ also, the matrix involved in solving this set of equations has rational entries. Indeed, the equation is equivalent to 
$$
(I_2\otimes X_1+X_1^T\otimes I_2) \Vec{X_{2}}=\Vec{I_2}.
$$
So the solution $X_{2}$ has entries in $\mathbb{Q}$. Solving this for $X_{2}$, we find
$$
X_{2} =\begin{bmatrix} \tfrac{1}{4} & -\tfrac{1}{8} \\ \tfrac{1}{2} & \tfrac{1}{4}\end{bmatrix},
$$
and hence
$$
Y=\begin{bmatrix} 1 & \frac{1}{2} & \tfrac{1}{4} & -\tfrac{1}{8} \\  -2 & 1 & \tfrac{1}{2} & \tfrac{1}{4} \\ 0 & 0 & 1 & \tfrac{1}{2} \\ 0 & 0 & -2 & 1\end{bmatrix}.
$$
One easily checks that indeed $Y^2=B_1$. 

As a next step, consider 
$$
B=\begin{bmatrix} A_1 & I_2 & 0 \\ 0 & A_1 & I_2 \\ 0 & 0 & A_1 \end{bmatrix}.
$$
To find a matrix $Z$ over $\mathbb{Q}$ for which $Z^2=B$, we can take
$$
Z=\begin{bmatrix} X_1 & X_{2} & X_{3} \\ 0& X_1 & X_{2} \\ 0 & 0 & X_1\end{bmatrix},
$$
with $X_1$ and $X_{2}$ as above. To find $X_{3}$ we have to solve
$$
X_1X_{3} +X_{2}^2+X_{3}X_1=0.
$$
Again this is a linear equation in the entries of $X_{3}$ which can be solved explicitly (and uniquely).
One obtains
$$
X_{3}=\tfrac{1}{32}X_1,
$$
so
$$
Z=\begin{bmatrix} 1 & \tfrac{1}{2} &\tfrac{1}{4} & -\tfrac{1}{8} & \tfrac{1}{32} & \tfrac{1}{64}\\
-2 & 1 &  \tfrac{1}{2} &\tfrac{1}{4} & -\tfrac{1}{16} & \tfrac{1}{32}\\
0 & 0 & 1 & \tfrac{1}{2} &\tfrac{1}{4} & -\tfrac{1}{8} \\
0 & 0 & -2 & 1 &  \tfrac{1}{2} &\tfrac{1}{4}\\
0 & 0 & 0 & 0 & 1 &  \tfrac{1}{2} \\
0 & 0 & 0 & 0 & -2 & 1 \end{bmatrix}.
$$

Now we are able to construct a rational solution $X$ to $X^2=A$. Indeed, if we recall that $B=T^{-1}AT$ it is clear that $X=TZT^{-1}$ satisfies $X^2=TZ^2T^{-1}=TBT^{-1}=A$. Computing $X$ explicity yields:
$$
X=\tfrac{1}{1024}
\begin{bmatrix} 
672 & 720 & -112 & 72 & -6 & 5 \\
-320 & 672 & 480 & -112 & 12 & -6 \\ 
384 & -320 & 960 & 480 & -40 & 12 \\ 
-768  & 384 & -896 & 960 & 336 & -40 \\
2560 & -768 & 2304 & -896 & 1440 & 336 \\
-21504 & 2560 & -16896 & 2304 & -4928 & 1440 
\end{bmatrix}.
$$
\hfill$\Box$

\medskip

\noindent
{\bf Example 2.} 
Consider the matrix
$$
A=\begin{bmatrix} -2 & 2 & 1 & 0 \\ -2 & -2 & 0 & 1 \\ 0 & 0 & -2 & 2 \\ 0 & 0 & -2 & -2 \end{bmatrix}
=\begin{bmatrix} A_1 & I_2 \\ 0 & A_1 \end{bmatrix}.
$$
We are interested in solving $X^3=A$, in particular, in finding a rational cube root of $A$, i.e., our polynomial is $p(\lambda)=\lambda^3$ with $\ell=3$ as in Section 8. Note that $A$ is in real Jordan canonical form, not in companion-Jordan canonical form. The transformation between the two is fairly easy. Let
$$
P=\begin{bmatrix} 1 & 0 \\ -2 & 2\end{bmatrix} \mbox{\ and \ } J=\begin{bmatrix} -2 & 2 \\ -2 & -2\end{bmatrix},
$$
then $PJP^{-1}=\begin{bmatrix} 0 & 1 \\ -8 & -4 \end{bmatrix}$. Hence, if we take $T=P\oplus P$, then $TAT^{-1}$ is in companion-Jordan canonical form.

\subsubsection*{Cube roots of the block upper left corner}

Put $A_1=\begin{bmatrix} 0 & 1 \\ -8 & -4 \end{bmatrix}$. As in the previous example, we first follow Drazin's approach to find a rational cube root of $A_1$. Let $\lambda_0=-2+2i$, which is one of the eigenvalues of $A_1$. The equation $\mu^3=\lambda_0$ has a solution in $\mathbb{Q}(i)$, namely $\mu=1+i$. Note that the other two cube roots are not in $\mathbb{Q}(i)$, but require $\mathbb{Q}(i,\sqrt{3})$. Indeed,
$$
\mu^3+(2-2i)=(\mu-(1+i))(\mu^2+(1+i)\mu+(1+i)^2)
$$
and so the other two cube roots are $-\frac{1}{2}(1+i)\pm \frac{1}{2}(1+i)i\sqrt{3}=
(1+i)(-\frac{1}{2}\pm\frac{1}{2}\sqrt{3}i)=(1+i)e^{k \pi i/3}$ with $k=2,4$ as expected. Note also that the characteristic polynomial of $A_1$ is reducible over $\mathbb{Q}(i)$ of course. 

The eigenvector of $A_1$ corresponding to the eigenvalue $\lambda_0$ is 
$$
\begin{bmatrix} 1 \\ -2 +2i\end{bmatrix} = \begin{bmatrix} 1 & 0 \\ -2 & 2 \end{bmatrix} \begin{bmatrix} 1 \\ i \end{bmatrix},
$$
so we have that $W=P$. To find a rational $X_1$ such that $X_1^3=A_1$ we need to solve $X_1W=M$, where $M$ is determined from $\mu \begin{bmatrix} 1 \\ -2 +2i\end{bmatrix} =\begin{bmatrix} 1+ i \\ -4 \end{bmatrix}=M\begin{bmatrix} 1 \\ i \end{bmatrix}$. Hence 
$
X_1=MW^{-1}=\begin{bmatrix}1 & 1 \\ -4 & 0 \end{bmatrix} \begin{bmatrix} 1 & 0 \\ 1 & \tfrac{1}{2}\end{bmatrix} =\begin{bmatrix} 2& \tfrac{1}{2} \\ -4 & 0 \end{bmatrix}.
$
It is easy to check that indeed $X_1^3=A_1$.

\subsubsection*{Rational cube roots of $A$}

As a next step we now find a rational cube root of $B=\begin{bmatrix} A_1 & I_2 \\ 0 & A_1\end{bmatrix}$. As before, we try 
$X=\begin{bmatrix} X_1 & X_{2} \\ 0 & X_1\end{bmatrix}$. Then
$$
X^3=\begin{bmatrix} X_1^2 & X_1X_{2}+X_{2}X_1\\ 0 & X_1^2\end{bmatrix} \begin{bmatrix} X_1 & X_{2} \\ 0 & X_1\end{bmatrix}=
\begin{bmatrix} X_1^3 & X_1^2X_{2} +X_1X_{2}X_1+X_{2}X_1^2
\\ 0 & X_1^3\end{bmatrix}.
$$
So the equation for the (1,2)-block is now
$$
X_1^2X_{2} +X_1X_{2}X_1+X_{2}X_1^2=I_2.
$$
Notice that although we now cannot make such a precise statement about solvability as in the square root case, we can say that if this is solvable at all, then the solution must have entries in the same field as where the entries of $X_1$ are located, as this is a linear equation in the entries of $X_{2}$. If we stack the columns of $X_{2}$ into $\Vec{X_{2}}$, the equation becomes 
$$
\left(I_2\otimes X_1^2 +X_1^T\otimes X_1+(X_1^2)^T\otimes I_2\right)\Vec{X_{2}}=\Vec{I_2}.
$$
Now we know necessary and sufficient conditions for an $m$th root to exist when we work over $\mathbb{C}$, and these are satisfied in this case. Then there is a solution, and by the argument above, the field over which a solution exists is dictated by the field over which we can find $X_1$. 

With $X_1$ as above we obtain
$$
I_2\otimes X_1^2 +X_1^T\otimes X_1+(X_1^2)^T\otimes I_2 =
\begin{bmatrix} 8 & 2 & -16 & -2 \\ -16 & 0 & 16 & -8 \\ 2 & \tfrac{1}{4} & 0 & 1 \\
-2 & 1 & -8 & -4\end{bmatrix},
$$
which is invertible. Hence $X_{2}$ is uniquely determined from the equation \eqref{eq:X2}.
Then 
$$
X_{2}= \tfrac{1}{12} \begin{bmatrix} -2  & -1 \\ 8 & 2\end{bmatrix}.
$$
Put $Y=\begin{bmatrix} X_1 & X_{2} \\ 0 & X_1\end{bmatrix}$, then $Y^3=B$, as one readily checks.

Finally, to find a rational cube root of $A$, note that $B=TAT^{-1}$, so with $X=T^{-1}YT$ we have 
$X^3=T^{-1}Y^3T=T^{-1}BT=A$. We find
$$
X=\begin{bmatrix} 1 & 1 & 0 & -\tfrac{1}{16} \\ 
-1 & 1 & \tfrac{1}{16} & 0 \\ 0 & 0 & 1 & 1 \\ 0 & 0 & -1 & 1\end{bmatrix}.
$$
\hfill $\Box$

\medskip

\noindent
{\bf Example 3.}
Consider the companion matrix 
$$
A=\begin{bmatrix} 0 & 1 & 0 & 0 & 0 & 0 \\ 0 & 0 & 1 & 0 & 0 & 0 \\ 0 & 0 & 0 & 1 & 0 & 0 \\
0 & 0 & 0 & 0 & 1 & 0 \\ 0 & 0 & 0 & 0 & 0 & 1\\ -4 & 0 & 0 & -4 & 0 & 0 \end{bmatrix},
$$
with characteristic polynomial $p_A(\lambda) = \lambda^6+4\lambda^3+4 = (\lambda^3+2)^2$. 
Find a rational solution $X$, such that $X^3 - 4X^2 +I = A$, that is, $p(X)=A$ where $p(\lambda)= \lambda^3-4\lambda^2+1$.

As in Example 1, we will construct a matrix $T$ such that $T^{-1}AT = B = \begin{bmatrix} A_1 & I_3 \\ 0 & A_1 \end{bmatrix}$, where $A_1 = \begin{bmatrix} 0 & 1 & 0 \\ 0 & 0 & 1 \\ -2 & 0 & 0 \end{bmatrix}$. We have $g(\lambda) = \lambda^3 + 2$, i.e., $d=2$ and $k=3$. Then $\lambda_0=-\sqrt[3]{2}$ is one of the roots of $p_A(\lambda)$ with associated eigenvector 
$$v_6(-\sqrt[3]{2}) = \begin{bmatrix} 1 \\ -\sqrt[3]{2} \\ \sqrt[3]{4} \\ -2 \\ 2\sqrt[3]{2} \\ -2\sqrt[3]{4} \end{bmatrix}.$$

Taking $\{1, -\sqrt[3]{2}, \sqrt[3]{4}\}$ as basis of $\mathbb{Q}(\sqrt[3]{2})$ over $\mathbb{Q}$, then  
$$v_6(-\sqrt[3]{2}) = \begin{bmatrix} 1 \\ -\sqrt[3]{2} \\ \sqrt[3]{4} \\ -2 \\ 2\sqrt[3]{2} \\ -2\sqrt[3]{4} \end{bmatrix}
=\begin{bmatrix} 1 & 0 & 0 \\ 0 & 1 & 0 \\ 0 & 0 & 1 \\ -2 & 0 & 0 \\ 0 & -2 & 0 \\ 0 & 0 & -2 \end{bmatrix}
\begin{bmatrix}1\\-\sqrt[3]{2}\\ \sqrt[3]{4}\end{bmatrix} = Wv_3(-\sqrt[3]{2}).$$
As before
$$S_1 = \begin{bmatrix} 0 & 0 & 0 & 0 & 0 & 0 \\ 1 & 0 & 0 & 0 & 0 & 0 \\ 0 & 2 & 0 & 0 & 0 & 0 \\ 0 & 0 & 3 & 0 & 0 & 0 \\ 0 & 0 & 0 & 4 & 0 & 0 \\0 & 0 & 0 & 0 & 5 & 0 \end{bmatrix},$$
so, $$
T = \begin{bmatrix} W & S_1\,W \end{bmatrix}=
\begin{bmatrix} 1 & 0 & 0 & 0 & 0 & 0 \\ 0 & 1 & 0 & 1 & 0 & 0 \\ 0 & 0 & 1 & 0 & 2 & 0 \\ -2 & 0 & 0 & 0 & 0 & 3 \\ 0 & -2 & 0 & -8 & 0 & 0 \\ 0 & 0 & -2 & 0 & -10 & 0 \end{bmatrix}.$$
An easy calculation verifies that indeed $T^{-1}AT = B$. We start by using Drazin's result to look for a solution $X_1$ such that $p(X_1) = A_1$ in the given field. If it exists, this result will be used to find a solution of $p(X)=A$. The characteristic polynomial $p_{A_1}(\lambda) = \lambda^3 + 2$ has as one of its roots $\lambda_0=-\sqrt[3]{2}$. Since $p(\lambda) = 1 - 4\lambda^2 +\lambda^3$, we are looking for a $\mu$ such that $p(\mu) = \lambda_0$. The desired $\mu$ is equal to $\mu = 1 + \sqrt[3]{2} + \sqrt[3]{4}$ in  $\mathbb{Q}(\sqrt[3]{2})$.

Following Drazin's approach, we take the eigenvector $v_3(\lambda_0)$ of $A_1$ corresponding to $\lambda_0$. Note that this gives $W=I_3$. Then we construct $M$ from $\mu v_3(\lambda_0)=Mv_3(\lambda_0)$:
$$
(1+\sqrt[3]{2}+\sqrt[3]{4})\begin{bmatrix} 1 \\ -\sqrt[3]{2} \\ \sqrt[3]{4}\end{bmatrix} =
\begin{bmatrix} 1+\sqrt[3]{2}+\sqrt[3]{4} \\ -2 -\sqrt[3]{2}- \sqrt[3]{4} \\ 2+2\sqrt[3]{2}+\sqrt[3]{4}\end{bmatrix}
= \begin{bmatrix} 1 & -1 & 1 \\ -2 & 1 & -1 \\ 2 & -2 & 1 \end{bmatrix} \begin{bmatrix} 1 \\ -\sqrt[3]{2} \\ \sqrt[3]{4}\end{bmatrix}.
$$
Then $X_1=MW^{-1} =M=\begin{bmatrix} 1 & -1 & 1 \\ -2 & 1 & -1 \\ 2 & -2 & 1 \end{bmatrix}$. It is easily checked that indeed $X_1^3-4X_1^2+I_3= A_1$.

Now to find $Y$ such that $p(Y)=\begin{bmatrix} A_1 & I_3 \\ 0 &A_1\end{bmatrix}$, we have to solve the following system of equations for $X_2$:
$$
X_2X_1^2+X_1X_2X_1+X_1^2X_2-4(X_1X_2+X_2X_1)=I_3,
$$
or equivalently,
$$
\left(I_3\otimes X_1^2+X_1^T\otimes X_1+(X_1^T)^2\otimes I_3 -4(X_1^T\otimes I_3+I_3\otimes X_1)\right) \Vec{X_2}=\Vec{I_3}.
$$
The solution is 
$$
X_2=\tfrac{1}{307}\begin{bmatrix} 41 & 26 & 9 \\ -18 & 41 & 26 \\ -52 & -18 & 41\end{bmatrix}.
$$
Taking $Y=\begin{bmatrix} X_1 & X_2 \\ 0 & X_1\end{bmatrix}$ yields a rational solution to 
$p(Y)=\begin{bmatrix} A_1 & I_3 \\ 0 &A_1\end{bmatrix}$.
Using the similarity to $A$ of the latter matrix, we can then find a rational solution $X$ to $p(X)=A$.
In fact, $X=TYT^{-1}$ will be a solution. We obtain after some calculation
$$
X=\tfrac{1}{3^3\cdot 307} 
\begin{bmatrix} 8451 & -11421 & 13581 & 81 & -1566 & 2646 \\
-10584 & 8451 & -11421 & 2997 &  81 & -1566\\
6264 & -10584 & 8451 & -5157 & 2997 &  81\\
-324 & 6264 & -10584 & 8127 &  -5157 & 2997 \\
-11988 & -324 & 6264 & -22572 & 8127 &  -5157  \\
20628 & -11988 & -324 &  26892 &  -22572 & 8127 
\end{bmatrix}.
$$
It is again checked, that $X^3-4X^2+I_6=A$. \hfill$\Box$

\medskip

\noindent
{\bf Example 4.}
Consider the case where $A$ is the $5\times 5$ standard Jordan block with eigenvalue $3$, and let the polynomial $p(\lambda)=\lambda^3-4\lambda^2 +5\lambda+1$. Then $p(1)=3, p(2)=3$, and $p^\prime(2)=1$, while $p^\prime(1)=0$. We see from this example that it can happen that for one solution $X_{1}$ of $p(X_{1})=C_{g}$ the condition \eqref{eq:detcondition} is satisfied, while for another solution it is not satisfied. 

To complete the example, consider the equations \eqref{eq:X2} and \eqref{eq:XjKron}. These are not solvable at all when we take the solution $X_1=1$, while for $X_1=2$ they become
\begin{align*}
p^\prime(2)X_2& =1, \mbox{\ so } X_2=1,\\
p^\prime(2)X_3& = -2 X_2^2, \mbox{\ so } X_3=-2.
\end{align*}
Next, we compute $X_4$ using the previously computed $X_1, X_2$ and $X_3$:
\begin{align*}
&p^\prime(2)X_4= -\sum_{m=2}^3 p_m \sum_{i=0}^{m-1}\sum_{\nu=0}^{m-2}
2^{m-2-\nu}\begin{bmatrix} 1 & -2 \end{bmatrix} \begin{bmatrix} 2& 1 \\ 0 & 2 \end{bmatrix}^\nu 
\begin{bmatrix}-2 \\ 1\end{bmatrix}
\\
=& -\left( \underbrace{-4 \cdot -4}_{m=2,i=0,\nu=0}+\underbrace{1 \cdot 2\cdot -4}_{m=3, i=0, \nu=0}
+\underbrace{1\cdot 1\cdot \begin{bmatrix} 1 & -2\end{bmatrix} \begin{bmatrix} -3 \\ 2 \end{bmatrix}}_{m=3, i=0, \nu=1} +\underbrace{1\cdot 2\cdot -4}_{m=3, i=1, \nu=0}\right) \\
=& -(16-8-7-8)=7.
\end{align*}
Finally, we compute $X_5$ using a similar computation
{\small
\begin{align*}
&p^\prime(2)X_5 = -\sum_{m=2}^3 p_m \sum_{i=0}^{m-1}\sum_{\nu=0}^{m-2}
2^{m-2-\nu}\begin{bmatrix} 1 & -2 & 7\end{bmatrix} \begin{bmatrix} 2& 1& -2  \\ 0 & 2& 1 \\
0 & 0 & 2 \end{bmatrix}^\nu 
\begin{bmatrix}7 \\-2 \\ 1\end{bmatrix} \\
=&
 -\left( \underbrace{-4 \cdot -18}_{m=2,i=0,\nu=0}+\underbrace{1 \cdot 2\cdot 18}_{m=3, i=0, \nu=0}+\underbrace{1 \cdot 2\cdot 18}_{m=3, i=1, \nu=0}
+\underbrace{1\cdot 1\cdot \begin{bmatrix} 1 & -2 & 7\end{bmatrix} \begin{bmatrix} 2& 1& -2  \\ 0 & 2& 1 \\
0 & 0 & 2 \end{bmatrix} 
\begin{bmatrix}7 \\-2 \\ 1\end{bmatrix}}_{m=3, i=1, \nu=1}\right) \\
=& - \begin{bmatrix} 1 & -2 & 7\end{bmatrix}\begin{bmatrix}10 \\-3 \\ 2\end{bmatrix}
=-30.
\end{align*}
}

\noindent
Hence, a rational solution to $X^3-4X^2+5X+I_5=J_5(3)$ is given by
$$
X=\begin{bmatrix} 2 & 1 & -2 & 7 & -30 \\
0 &  2 & 1 & -2 & 7 \\
0 & 0 &  2 & 1 & -2 \\
0 & 0 & 0 & 2 & 1 \\ 
0 & 0 & 0 & 0 & 2 \end{bmatrix}.
$$
\hfill$\Box$

\paragraph{Acknowledgements:}
This work is based on research supported in part by the DSI-NRF Centre of Excellence in Mathematical and Statistical Sciences (CoE-MaSS) and in part by the National Research Foundation of South Africa (Grant Numbers 145688 and 2022-012-ALG-ILAS).
Opinions expressed and conclusions arrived at are those of the authors and are not necessarily to be attributed to the CoE-MaSS.

\end{document}